\documentclass[11pt]{article}
\usepackage{amsmath}
\usepackage{amsthm}
\usepackage{amssymb}
\usepackage{mathrsfs}
\usepackage{amssymb,amsmath}
\usepackage{leftidx}
\usepackage{graphicx}
\usepackage{amsfonts}
\usepackage{graphicx, color}
\usepackage{cite}
\usepackage{enumerate}
\usepackage{enumitem}
\usepackage{subfigure}
\usepackage{booktabs}
\usepackage{algpseudocode,algorithm,algorithmicx}

\newtheorem{theorem}{Theorem}[section]
\newtheorem{lemma}{Lemma}[section]

\newtheorem{example}{Example}
\newtheorem{remark}{Remark}

\newcommand{\zd}{\,\mathrm{d}}

\newcommand{\diff}{\triangledown_{\tau}}
\newcommand{\abs}[1]{\left|#1\right|}
\newcommand{\absb}[1]{\big|#1\big|}

\newcommand{\abst}[1]{|#1|}
\newcommand{\bra}[1]{\left(#1\right)}
\newcommand{\brab}[1]{\big(#1\big)}
\newcommand{\braB}[1]{\Big(#1\Big)}
\newcommand{\brat}[1]{(#1)}
\newcommand{\kbra}[1]{\left[#1\right]}
\newcommand{\kbrab}[1]{\big[#1\big]}
\newcommand{\kbraB}[1]{\Big[#1\Big]}

\newcommand{\myinner}[1]{\left\langle#1\right\rangle}
\newcommand{\myinnerb}[1]{\big\langle#1\big\rangle}
\newcommand{\myinnerB}[1]{\Big\langle#1\Big\rangle}
\newcommand{\mynorm}[1]{\left\|#1\right\|}
\newcommand{\mynormb}[1]{\big\|#1\big\|}

\newcommand{\timenorm}[1]{\absb{\!\absb{\!\absb{#1}\!}\!}}
\newcommand{\myvec}[1]{\boldsymbol{#1}}

\begin{document}
\title{Mesh-robustness of an energy stable BDF2 scheme with  \\variable steps
for the Cahn-Hilliard model}
\author{
Hong-lin Liao\thanks{ORCID 0000-0003-0777-6832; Department of Mathematics,
Nanjing University of Aeronautics and Astronautics,
Nanjing 211106, China; Key Laboratory of Mathematical Modelling
and High Performance Computing of Air Vehicles (NUAA), MIIT, Nanjing 211106, China.
Hong-lin Liao (liaohl@nuaa.edu.cn and liaohl@csrc.ac.cn)
is supported by a grant 12071216 from
National Natural Science Foundation of China.}
\quad
Bingquan Ji\thanks{Department of Mathematics, Nanjing University of Aeronautics and Astronautics,
211101, P. R. China. Bingquan Ji (jibingquanm@163.com).}
\quad
Lin Wang\thanks{Beijing Computational Science Research Center, Beijing 100193, P. R. China. E-mail: wanglin@csrc.ac.cn.}
\quad Zhimin Zhang\thanks{Beijing Computational Science Research Center, Beijing 100193, P. R. China;
and Department of Mathematics, Wayne State University, Detroit, MI 48202, USA.
E-mails: zmzhang@csrc.ac.cn  and ag7761@wayne.edu.
This author is supported in part by the NSFC grant 11871092 and NSAF grant U1930402.}
}
\date{\today}
\maketitle
\normalsize

\begin{abstract}
  The two-step backward differential formula (BDF2) with unequal time-steps
  is applied to construct an energy stable convex-splitting scheme for
  the Cahn-Hilliard model.   We focus on the numerical influences of time-step variations
  by using the recent theoretical framework
  with the discrete orthogonal convolution kernels.
  Some novel discrete convolution embedding inequalities with respect to
  the orthogonal convolution kernels are developed such that
  a concise $L^2$ norm error estimate is established  at the first time
  under an updated step-ratio restriction   $0 <r_k:=\tau_k/\tau_{k-1}\leq r_{\mathrm{user}}$,
  where $r_{\mathrm{user}}$ can be chosen by the user   such that $r_{\mathrm{user}}<4.864$.
  The stabilized convex-splitting BDF2 scheme is shown to be mesh-robustly convergent
  in the sense that the convergence constant (prefactor)
  in the error estimate is independent of the adjoint time-step ratios.
  The suggested method is proved to preserve a modified energy dissipation law
  at the discrete levels if $0<r_k\le r_{\mathrm{user}}$,
  such that it is mesh-robustly stable in an energy norm.
  On the basis of ample tests on random time meshes,
  a useful adaptive time-stepping strategy is applied to efficiently capture
  the multi-scale behaviors and to accelerate
  the long-time simulation approaching the steady state.
\end{abstract}

\section{Introduction}\setcounter{equation}{0}

The Cahn-Hilliard (CH) model is an efficient approach to
describe the coarsening dynamics of a binary alloy system \cite{Cahn1958Free}
and has been applied in other fields including image inpainting \cite{Bertozzi2007Inpainting}
and tumor growth \cite{Cristini2009Nonlinear}.
Consider a free energy functional of Ginzburg--Landau type,
\begin{align}\label{cont:free energy}
E[\Phi] = \int_{\Omega}\kbraB{\frac{\epsilon^2}{2}\abst{\nabla\Phi}^2+F(\Phi)}\zd\myvec{x}\quad\text{with}\quad F(\Phi):=\frac14(\Phi^2-1)^2
\end{align}
where $\myvec{x}\in\Omega\subseteq\mathbb{R}^2$ and $0<\epsilon<1$ is a bounded parameter that is proportional to the interface width.
Then the Cahn-Hilliard equation would be given by the $H^{-1}$ gradient flow
associated with the free energy functional $E[\Phi]$,
\begin{align}\label{cont: Problem-CH}
\partial_t \Phi=\kappa\Delta\mu\quad\text{with}\quad
\mu:=\tfrac{\delta E}{\delta \Phi}
=F'(\Phi)-\epsilon^2\Delta\Phi,
\end{align}
where the parameter $\kappa$ is the mobility related to the characteristic
relaxation time of system and $\mu$ is the chemical potential.
Assume that $\Phi$ is periodic over the domain $\Omega$.
By applying the integration by parts, one can find the volume conservation,
\begin{align}\label{cont:volume conservation}
\brab{\Phi(t),1}=\brab{\Phi(t_0),1},
\end{align}
and the following energy dissipation law,
\begin{align}\label{cont:energy dissipation}
\frac{\zd{E}}{\zd{t}}=\brab{\tfrac{\delta E}{\delta \Phi},\partial_t \Phi}
=\kappa\bra{\mu,\Delta\mu}
=-\kappa\mynorm{\nabla\mu}^2 \le 0,
\end{align}
where $\bra{u,v}:=\int_{\Omega}uv\zd{\myvec{x}}$, and the associated
$L^{2}$ norm $\mynorm{v}=\sqrt{\bra{v,v}}$ for all $u,v\in{L}^{2}(\Omega)$.

The main aim of this paper is to present a rigorous stability and convergence analysis
of the BDF2  method with variable time-steps for simulating the CH model \eqref{cont: Problem-CH}.
Consider the nonuniform time levels $0=t_{0}<t_{1}<\cdots<t_N=T$
with the time-step sizes $\tau_{k}:=t_{k}-t_{k-1}$ for $1\le k \le N$,
and denote the maximum time-step size $\tau:=\max_{1\le k\le N}\tau_k$.
Let the adjoint time-step ratio $r_k:=\tau_k/\tau_{k-1}$ for $2\le k\le N$.
Our analysis will focus on the influence of non-uniform time grids (with the associated time-step ratios)
on the numerical solution by carefully evaluating  the stability and convergence.

This is motivated by the following facts:
\begin{itemize}
\item The BDF2 method is A-stable and L-stable such that it would be more suitable than Crank-Nicolson type schemes
   for solving the stiff dissipative problems,
  see e.g. \cite{ChengFengWangWise:2019Fourth,ChengWangWise:2019An}.
  \item The nonuniform grid and adaptive time-stepping techniques  \cite{GomezHughes:2011Provably,LiaoJiZhang:2020pfc,LiaoSongTangZhou:2020mbe,LiaoTangZhou:2020maximumAC,QiaoZhangTang:2011}
  are powerful in capturing the multi-scale behaviors
  and accelerating the long-time simulations of phase field models including the CH model.
  \item The convergence theory of variable-steps BDF2 scheme remains incomplete
  for nonlinear parabolic equations. Actually, the required step-ratio constraint
  for the $L^2$ norm stability are severer than the classical zero-stability condition
  $r_k<1+\sqrt{2}$, given by Grigorieff \cite{Grigorieff:1983}.
  Always, they contain some undesirable pre-factors $C_r\exp(C_r\Gamma_n)$ or $C_r\exp(C_rt_n)$, see  e.g.
  \cite{Becker:1998,ChenWangYanZhang:2019,Emmrich:2005,WangChenFang:2019},
  where $\Gamma_n$ may be unbounded when certain time-step variations appear
  and $C_r$ may be infinity as the step-ratios approach the zero-stability limit $1+\sqrt{2}$.
\end{itemize}
In recent works \cite{LiaoJiZhang:2020pfc,LiaoSongTangZhou:2020mbe,LiaoZhang:2020linear},
a novel technique with \emph{discrete orthogonal convolution} (DOC) kernels
was suggested to verify that, if
$0< r_k < \bra{3+\sqrt{17}}/2 \approx 3.561$,
the BDF2 scheme is computationally robust with respect to
the time-step variations for linear diffusions \cite{LiaoZhang:2020linear},
the phase field crystal model \cite{LiaoJiZhang:2020pfc}
and the molecular beam epitaxial model without slope selection \cite{LiaoSongTangZhou:2020mbe}.

Nonetheless, due to the lack of some convolution embedding inequalities
with respect to the DOC kernels, the techniques
in \cite{LiaoJiZhang:2020pfc,LiaoSongTangZhou:2020mbe,LiaoZhang:2020linear}
are inadequate to handle more general nonlinear problems such as
the underlying nonlinear CH model (and Allen-Chan model).
The main aim of this paper is to fill this gap by establishing some discrete convolution
embedding inequalities with respect to the DOC kernels.
Also, the recent analysis in \cite[Lemma A.1]{LiaoJiZhang:2020pfc}
with a step-scaled matrix motivates us
to update the previous zero-stability restriction in \cite{LiaoZhang:2020linear} as follows,
\begin{enumerate}[itemindent=1em]
\item[$\mathbf{S0}$.]  $0< r_k \le r_{\mathrm{user}}(<4.864)$   \ for \ $2\le k \le N$,
\end{enumerate}
where the value of $r_{\mathrm{user}}$ can be chosen in adaptive time-stepping computations
by the user such that $r_{\mathrm{user}}<4.864$, such as $r_{\mathrm{user}}=2,3$ or 4 for practical choices.
Under the step-ratio constraint $\mathbf{S0}$,
we will present an $L^2$ norm error estimate with an improved prefactor, see Theorem \ref{thm: L2 Convergence-CH},
$$C_{\phi}\exp\brab{c_{\epsilon} t_{n-1}}.$$
Here and hereafter, any subscripted $C$, such as $C_u$ and $C_\phi$, denotes a generic positive constant, not necessarily
 the same at different occurrences; while,
 any subscripted $c$, such as $c_{\epsilon}, c_\Omega,c_p, c_z$ and so on,
denotes a fixed constant. The appeared constants may be dependent on the given data
(typically, the interface width parameter $\epsilon$)
 and the solution but are always independent of the spatial lengths, the time $t_n$, the step sizes $\tau_n$ and the step ratios $r_n$.
It is interesting to emphasize that, under the step-ratio constraint $\mathbf{S0}$,
the involved constants are bounded even when the step-ratios $r_n$
approach  $r_{\mathrm{user}}$ such that the BDF2 scheme is mesh-robustly convergent.

To the best of our knowledge, this is the first time such an optimal $L^2$ norm error estimate of variable-steps BDF2 method
 is established for the Cahn-Hiliard (and Allen-Cahn) type models.
As a closely related work, the BDF2 scheme for the Allen-Chan equation was also investigated
in \cite{LiaoTangZhou:2020maximumAC} by using the discrete complementary convolution kernels.
The BDF2 scheme was proved to preserve the maximum bound principle if the step-ratios satisfy
the classical zero-stability condition $r_k<1+\sqrt{2}$. The maximum norm error estimate with a prefactor $\frac{1}{1-\eta}\exp(\frac{t_n}{1-\eta})$ was obtained,
where the parameter $\eta\rightarrow1$ as $\max r_k \rightarrow1+\sqrt{2}$.
It is to mention that,
under the constraint $\mathbf{S0}$,
one can follow the present analysis to
obtain a new $L^2$ norm error estimate that is robustly stable to
the variations of time-steps.

Given a grid function $\{v^k\}_{k=0}^N$,
put $\diff v^{k}:=v^{k}-v^{k-1}$, $\partial_{\tau}v^{k}:=\diff v^{k}/\tau_k$ for $k\geq{1}$.
Taking $v^n=v(t_n)$, we view the variable-steps BDF2 formula as a discrete convolution summation
\begin{align}\label{def: BDF2-Formula}
D_2v^n:=\sum_{k=1}^nb_{n-k}^{(n)}\diff v^k
\quad \text{for $n\ge2$},
\end{align}
where the discrete convolution kernels $b_{n-k}^{(n)}$ are defined for $n\ge2$,
\begin{align}\label{def: BDF2-kernels}
b_{0}^{(n)}:=\frac{1+2r_n}{\tau_n(1+r_n)},\quad
b_{1}^{(n)}:=-\frac{r_n^2}{\tau_n(1+r_n)}\quad \text{and} \quad
b_{j}^{(n)}:=0\quad \mathrm{for}\quad j\ge2.
\end{align}
Without losing the generality,
assume that an accurate solution $\phi^1$ is available.
We consider the stability and convergence of
the convex-splitting BDF2 scheme
for solving the CH equation \eqref{cont: Problem-CH}
subject to the periodic boundary conditions:
\begin{align}\label{scheme: BDF2 convex CH}
D_2\phi^n=\kappa\Delta_h\mu^n\quad\text{with}
\quad\mu^n:=\brab{\phi^n}^3-\hat{\phi}^n-\bra{\epsilon^2
+A\tau^2}\Delta_h\phi^n
\quad\text{for $2\le n\le N$},
\end{align}
where $\hat{\phi}^n:=\bra{1+r_n}\phi^{n-1}-r_n\phi^{n-2}$
and the stabilized parameter $A>0$.
The spatial operators are approximated by the Fourier pseudo-spectral method,
as described in the next section.

The unique solvability of the convex-splitting scheme \eqref{scheme: BDF2 convex CH}
is established in Theorem \ref{thm: convexity solvability CH}
by using the fact that the solution of nonlinear scheme \eqref{scheme: BDF2 convex CH}
is equivalent to the minimization of a convex functional.
Lemma \ref{lem: BDF2-Kernels-Positive} shows that the BDF2 convolution
kernels $b_{n-k}^{(n)}$ are positive definite provided the adjacent time-step
rations $r_k$ satisfy $\mathbf{S0}$.
Theorem \ref{thm: energy-decay-law CH} shows that
the convex-splitting BDF2 method \eqref{scheme: BDF2 convex CH}
has a modified energy dissipation law at the discrete levels
for a properly large parameter $A$,
see Remark \ref{remark: comments on stabilized A condition}.

We are to emphasize that the solution estimates in section 2
are based on the original form \eqref{scheme: BDF2 convex CH}, but in the subsequent $L^2$ norm error analysis
we will use an equivalent convolution form with a class
of discrete orthogonal convolution (DOC)
kernels. The DOC kernels $\{\theta_{n-k}^{(n)}\}_{k=2}^n$ are defined by
(this definition is slightly different from those
in \cite{LiaoJiZhang:2020pfc,LiaoSongTangZhou:2020mbe,LiaoZhang:2020linear}
since we do not introduce the discrete kernel $b_{0}^{(1)}$
for the first-level solver)
\begin{align}\label{def: DOC-Kernels}
\theta_{0}^{(n)}:=\frac{1}{b_{0}^{(n)}}
\;\; \text{for $n\ge2$}\quad \mathrm{and} \quad
\theta_{n-k}^{(n)}:=-\frac{1}{b_{0}^{(k)}}
\sum_{j=k+1}^n\theta_{n-j}^{(n)}b_{j-k}^{(j)}
\quad \text{for $n\ge k+1\ge3$}.
\end{align}
One has the following  discrete orthogonal identity
   \begin{align}\label{eq: orthogonal identity}
 \sum_{j=k}^n\theta_{n-j}^{(n)}b_{j-k}^{(j)}\equiv \delta_{nk}\quad\text{for $2\le k\le n$,}
   \end{align}
where $\delta_{nk}$ is the Kronecker delta symbol. By exchanging the summation order
and using the identity \eqref{eq: orthogonal identity}, it is not difficult to check that
  \begin{align}\label{eq: orthogonal equality for BDF2}
  \sum_{j=2}^n\theta_{n-j}^{(n)}D_2v^j=&\,\sum_{j=2}^{n}\theta_{n-j}^{(n)}b_{j-1}^{(j)}\diff v^1
+\sum_{j=2}^{n}\theta_{n-j}^{(n)}
\sum_{\ell=2}^{j}b_{j-\ell}^{(j)}\diff v^\ell\nonumber\\
=&\,\theta_{n-2}^{(n)}b_{1}^{(2)}\diff v^1+\diff v^n\quad\text{for $n\ge2$.}
   \end{align}
Acting the DOC kernels $\theta_{m-n}^{(m)}$ on the first equation in \eqref{scheme: BDF2 convex CH}
and summing $n$ from $n=2$ to $m$, we apply \eqref{eq: orthogonal equality for BDF2}
to find the equivalent convolution form (replacing $m$ by $n$)
\begin{align}\label{scheme: DOC form BDF2 convex CH}
\diff\phi^n=-\theta_{n-2}^{(n)}b_{1}^{(2)}\diff \phi^1+\kappa\sum_{j=2}^n\theta_{n-j}^{(n)}\Delta_h\mu^j
\quad\text{for $2\le n\le N$.}
\end{align}
Note that, by following the proof of \cite[Lemma 2.1]{LiaoTangZhou:2020doc}, we have
 \begin{align}\label{eq: mutual orthogonal identity}
 \sum_{j=k}^mb_{m-j}^{(m)}\theta_{j-k}^{(j)}\equiv \delta_{mk}\quad\text{for $2\le k\le m$.}
   \end{align}
With the help of this mutually orthogonal identity,
one can recover the original form \eqref{scheme: BDF2 convex CH} by acting the BDF2 kernels $b_{m-n}^{(m)}$
on the new formulation \eqref{scheme: DOC form BDF2 convex CH}. In this sense, the DOC kernels define a \textit{reversible discrete transform}
between \eqref{scheme: BDF2 convex CH} and the convolution form \eqref{scheme: DOC form BDF2 convex CH}.

To perform the $L^2$ norm error analysis, section 3 presents some properties of
the DOC kernels $\theta_{n-k}^{(n)}$ and  some new convolution embedding inequalities with respect to the DOC kernels,
see Lemmas 3.1--3.9.
By making use of the $H^1$ norm solution bound obtained in Lemma \ref{lem: bound-BDF2 Solution CH},
we establish an optimal $L^2$ norm error estimate in section 4.
Numerical tests and comparisons are presented in section 5 to validate the
accuracy and effectiveness of the BDF2 method \eqref{scheme: BDF2 convex CH},
especially when coupled with an adaptive stepping strategy.

\section{Solvability and energy dissipation law}
\setcounter{equation}{0}

We use the same spatial notations in \cite{LiaoJiZhang:2020pfc}.
Set the space domain $\Omega=(0,L)^2$
and consider the uniform length $h_x=h_y=h:=L/M$ in each direction
for an even positive integer $M$.
Let $\Omega_{h}:=\big\{\myvec{x}_{h}=(ih,jh)\,|\,1\le i,j \le M\big\}$
and put
$\bar{\Omega}_{h}:=\Omega_{h}\cup\partial{\Omega}$.
Denote the space of $L$-periodic grid functions
$\mathbb{V}_{h}:=\{v\,|\,v=\bra{v_h}\; \text{is $L$-periodic for}\; \myvec{x}_h\in\bar{\Omega}_h\}.$

For a periodic function $v(\myvec{x})$ on $\bar{\Omega}$,
let $P_M:L^2(\Omega)\rightarrow \mathscr{F}_M$
be the standard $L^2$ projection operator onto the space $\mathscr{F}_M$,
consisting of all trigonometric polynomials of degree up to $M/2$,
and $I_M:L^2(\Omega)\rightarrow \mathscr{F}_M$
be the trigonometric interpolation operator \cite{ShenTangWang:2011Spectral},
i.e.,
\[
\bra{P_Mv}(\myvec{x})=\sum_{\ell,m =- M/2}^{M/2-1}
\widehat{v}_{\ell,m}e_{\ell,m}(\myvec{x}),\quad
\bra{I_Mv}(\myvec{x})=\sum_{\ell,m=- M/2}^{M/2-1}
\widetilde{v}_{\ell,m}e_{\ell,m}(\myvec{x}),
\]
where the complex exponential basis function
$e_{\ell,m}(\myvec{x}):=e^{\mathrm{i}\nu\bra{\ell x+my}}$ with $\nu=2\pi/L$.
The coefficients $\widehat{v}_{\ell,m}$
refer to the standard Fourier coefficients of function $v(\myvec{x})$,
and the
pseudo-spectral coefficients $\widetilde{v}_{\ell,m}$ are determined such that $\bra{I_Mv}(\myvec{x}_h)=v_h$.

The Fourier pseudo-spectral first and second order derivatives of $v_h$ are given by
\[
\mathcal{D}_xv_h:=\sum_{\ell,m= -M/2}^{M/2-1}
\bra{\mathrm{i}\nu\ell}\widetilde{v}_{\ell,m}
e_{\ell,m}(\myvec{x}_h),\quad
\mathcal{D}_x^2v_h:=\sum_{\ell,m = -M/2}^{M/2-1}
\bra{\mathrm{i}\nu\ell}^2\widetilde{v}_{\ell,m}
e_{\ell,m}(\myvec{x}_h).
\]
The differentiation  operators $\mathcal{D}_y$
and $\mathcal{D}_y^2$ can be defined in the similar fashion.
In turn, we can define the discrete gradient and Laplacian
in the point-wise sense, respectively, by
\[
\nabla_hv_h := \left(\mathcal{D}_xv_h,\mathcal{D}_yv_h\right)^T\quad\text{and}\quad
\Delta_hv_h :=\bra{\mathcal{D}_x^2+\mathcal{D}_y^2}v_h.
\]

For any grid functions $v,w\in\mathbb{V}_{h}$,
define the discrete inner product
$\myinner{v,w}:=h^2\sum_{\myvec{x}_h\in\Omega_{h}}v_h w_h$,
and the associated $L^{2}$ norm $\mynorm{v}:=\mynorm{v}_{l^2}=\sqrt{\myinner{v,v}}$.
Also, we will use the discrete $l^q$ norm $\mynorm{v}_{l^q}:=\sqrt[q]{h^2\sum_{\myvec{x}_h\in\Omega_{h}}|v_h|^q}$
and the $H^1$ seminorm $\mynormb{\nabla_hv}:=\sqrt{h^2\sum_{\myvec{x}_h\in\Omega_{h}}|\nabla_hv_h|^2}$.
It is easy to check the discrete Green's formulas,
$\myinner{-\Delta_hv,w}=\myinner{\nabla_hv,\nabla_hw}$ and
$\myinner{\Delta_h^2v,w}=\myinner{\Delta_hv,\Delta_hw}$,
see \cite{ShenTangWang:2011Spectral,ChenWangWang:2014,ChengWangWise:2019An} for more details.
Also we have the following discrete embedding inequality simulating the Sobolev embedding $H^1(\Omega)\hookrightarrow L^6(\Omega)$,
\begin{align}\label{ieq: H1 embedding L6}
\mynormb{v}_{l^6}\le c_\Omega\brab{\mynormb{v}+\mynormb{\nabla_hv}}\quad \text{for any $v\in\mathbb{V}_{h}$.}
\end{align}

For the underlying volume-conservative problem, it is also to define a mean-zero function space
$\mathbb{\mathring V}_{h}:=\big\{v\in\mathbb{V}_{h}\,|\, \myinner{v,1}=0\big\}\subset\mathbb{V}_{h}.$
As usual, following the arguments in \cite{ChengWangWiseYue:2016Weakly,ChengWangWise:2019An}, one can introduce an discrete version of
inverse Laplacian operator $\bra{-\Delta_h}^{-\gamma}$ as follows.
For a grid function $v\in\mathbb{\mathring V}_{h}$, define
\[
\bra{-\Delta_h}^{-\gamma}v_h
:=\sum_{\mbox{\tiny$\begin{array}{c}
\ell,m=-M/2\\
\bra{\ell,m}\neq \mathbf{0}
\end{array}$}}^{M/2-1}
\bra{\nu^2\bra{\ell^2+m^2}}^{-\gamma}\widetilde{v}_{\ell,m}
e_{\ell,m}(\myvec{x}_h),
\]
and an $H^{-1}$ inner product  $\myinner{v,w}_{-1}
:=\myinnerb{\bra{-\Delta_h}^{-1}v,w}.$
The associated $H^{-1}$ norm $\mynorm{\cdot}_{-1}$ can be defined by $\mynorm{v}_{-1}:=\sqrt{\myinner{v,v}_{-1}}\,.$
We have the following Poincar\'{e} type inequality with the usual Poincar\'{e} constant $c_p$, $\mynormb{v}_{-1}\le c_p\mynormb{v}$,
and the generalized H\"{o}lder inequality,
\begin{align}\label{ieq: generalized Holder H-1}
\mynormb{v}^2\le \mynormb{\nabla_hv}\mynormb{v}_{-1}\quad \text{for any $v\in\mathbb{\mathring V}_{h}$.}
\end{align}
Also the discrete embedding inequality \eqref{ieq: H1 embedding L6} can be simplified as ($c_z:=c_\Omega+c_\Omega c_p$)
\begin{align}\label{ieq: mean-zero H1 embedding L6}
\mynormb{v}_{l^6}\le c_z\mynormb{\nabla_hv}\quad \text{for any $v\in\mathbb{\mathring V}_{h}$.}
\end{align}

\subsection{Unique solvability}

Let $E[\phi^k]$ be the discrete version of free energy functional
\eqref{cont:free energy}, given by
\begin{align}\label{def: discrete free energy}
E[\phi^k]
:=\frac{\epsilon^2}{2}\mynormb{\nabla_h\phi^k}^2+
\myinnerb{F(\phi^k),1} \quad\text{for $k\ge 1$.}
\end{align}
To focus on the numerical analysis of the BDF2 solution, it is to assume that
\begin{enumerate}[itemindent=1em]
\item[$\mathbf{A1}$.]
 A starting scheme is properly chosen to compute the first-level solution $\phi^1$
 such that it preserves the volume,
  $\myinnerb{\phi^1,1}=\myinnerb{\phi^0,1}=\myinnerb{P_M\Phi^0,1}$,
 and also preserves certain (maybe, modified) energy dissipation law.
 There exists a positive constant $c_0$, depended on the domain $\Omega$,
  the mobility $\kappa$,
  the interface  parameter $\epsilon$ and the initial value $\phi^0$, such that
 $$E[\phi^1]+\frac{\tau_{2}}{2\kappa}\mynormb{\partial_{\tau} \phi^1}_{-1}^2
+\frac{\tau_1\tau_{2}}{2}\mynormb{\partial_{\tau}\phi^1}^2
+\frac{A\tau^2}{2}\mynormb{\nabla_h\phi^1}^2
 \le c_0.$$
 \end{enumerate}

 \begin{remark}\label{remark: first-level TR-BDF2}
Assumption $\mathbf{A1}$ can be satisfied by many of first-level solvers.
The BDF1 scheme would be suited for computing a second-order solution $\phi^1$; however,
a very small initial step $\tau_1$ would not be suggested here since it arrives at a large step-ratio $r_2$
and eventually affects the accuracy of solution in the whole simulation, see numerical results in \cite{Nishikawa:2019}.

The Crank-Nicolson scheme at the first time-level can generate a second-order difference quotient $\partial_\tau\phi^1$;
but a very small initial step $\tau_1$ would not be suggested either because it would be prone to generate nonphysical oscillations.
To control possibly initial oscillations, we suggest a special step-ratio $r_2=\sqrt{2}/2$ in the implementation of
our scheme \eqref{scheme: BDF2 convex CH}. Actually, by taking $\phi^{\gamma}:=\phi^1$, $\phi^1:=\phi^2$, $\tau_{*}:=\tau_1+\tau_2$
and $\gamma:=\tau_1/\tau_{*}$ with $r_2=1/\gamma-1$ , the first two steps of \eqref{scheme: BDF2 convex CH} are
equivalent to the following TR-BDF2 method
\begin{align*}
\frac{\phi^{\gamma}-\phi^0}{\gamma\tau_{*}}
=&\,\frac\kappa2\Delta_h\mu^{\gamma}+\frac\kappa2\Delta_h\mu^{0},\quad
\frac{2-\gamma}{(1-\gamma)\tau_{*}}\phi^{1}
-\frac{1}{\gamma(1-\gamma)\tau_{*}}\phi^{\gamma}
+\frac{1-\gamma}{\gamma\tau_{*}}\phi^{0}=\kappa\Delta_h\mu^{1},
\end{align*}
 which was shown to be L-stable for $\gamma=2-\sqrt{2}$, see \cite{HoseaShampine:1996,TumoloBonaventura:2015}.

\end{remark}

Under the assumption $\mathbf{A1}$,
the solution $\phi^n$ of the BDF2 scheme \eqref{scheme: BDF2 convex CH} preserves the volume, $\myinnerb{\phi^n,1}=\myinnerb{\phi^0,1}$
for $n\ge2$. Actually, taking the inner product of \eqref{scheme: BDF2 convex CH} by 1 and applying the discrete Green's formulas,
one can check that $\myinnerb{D_2\phi^j,1}=0$ for $j\ge2$.
 Multiplying both sides of this equality by the DOC kernels
 $\theta_{n-j}^{(n)}$ and summing the index $j$ from $j=2$ to $n$, we get
$$\sum_{j=2}^n\theta_{n-j}^{(n)}\myinnerb{D_2\phi^j,1}=0\quad\text{for $n\ge2$}.$$
It leads to $\myinnerb{\diff \phi^n,1}=0$ directly by taking $v^j=\phi^j$ in the equality \eqref{eq: orthogonal equality for BDF2}.
Simple induction yields the volume conversation law, $\myinnerb{\phi^n,1}=\myinnerb{\phi^{n-1},1}=\cdots=\myinnerb{\phi^{0},1}$
for $n\ge1$.

\begin{theorem}\label{thm: convexity solvability CH}
If $\mathbf{A1}$ holds, the convex-splitting BDF2 scheme
\eqref{scheme: BDF2 convex CH} is uniquely solvable.
\end{theorem}
\begin{proof}
For any fixed time-level indexes $n\ge2$,
consider the following energy functional $G$ on the space
$\mathbb{V}_{h}^{*}:=\big\{z\in\mathbb{V}_{h}\,|\, \myinnerb{z,1}=\myinnerb{\phi^{n-1},1}\big\},$
\begin{align}\label{eq: discrete energy functional}
G[z]:=&\frac{b_0^{(n)}}{2}\mynormb{z-\phi^{n-1}}_{-1}^2
+b_1^{(n)}\myinnerb{\diff \phi^{n-1},z-\phi^{n-1}}_{-1}\nonumber\\
&\quad+\frac{\kappa}{2}\brab{\epsilon^2+A\tau^2}\mynormb{\nabla_hz}^2
+\kappa\myinnerb{z^3/4-\hat{\phi}^n,z}.
\end{align}
It is easily to verity
the functional $G$ is strictly convex since,
for any $\lambda\in \mathbb{R}$ and any $\psi\in \mathbb{\mathring V}_{h}$,
\begin{align*}
\frac{\zd^2G}{\zd\lambda^2}[z+\lambda\psi]\Big|_{\lambda=0}
=&\,b_0^{(n)}\mynormb{\psi}_{-1}^2
+\kappa\brab{\epsilon^2+A\tau^2}\mynormb{\nabla_h\psi}^2
+3\kappa\mynormb{z\psi}^2>0.
\end{align*}
Thus the functional $G$ has a unique minimizer,
denoted by $\phi^n$, if and only if it solves the equation
\begin{align*}
0=\frac{\zd G}{\zd\lambda}[z+\lambda\psi]\Big|_{\lambda=0}=&\,
\myinnerb{b_0^{(n)}\brat{z-\phi^{n-1}}+b_1^{(n)}\diff \phi^{n-1},\psi}_{-1}
+\kappa\myinnerb{z^3-\hat{\phi}^n-\brab{\epsilon^2+A\tau^2}\Delta_hz,\psi}\\
=&\,\myinnerB{b_0^{(n)}\brat{z-\phi^{n-1}}+b_1^{(n)}\diff \phi^{n-1}
-\kappa\Delta_h\kbrab{z^3-\hat{\phi}^n-\brab{\epsilon^2+A\tau^2}\Delta_hz},\psi}_{-1}.
\end{align*}
This equation holds for any $\psi\in \mathbb{\mathring V}_{h}$ if and only
if the unique minimizer $\phi^n\in\mathbb{V}_{h}^{*}$ solves
\begin{align*}
b_0^{(n)}\brat{\phi^n-\phi^{n-1}}+b_1^{(n)}\diff \phi^{n-1}
-\kappa\Delta_h\kbrab{(\phi^n)^3-\hat{\phi}^n
-\brab{\epsilon^2+A\tau^2}\Delta_h\phi^n}=0,
\end{align*}
which is just the convex-splitting BDF2 scheme \eqref{scheme: BDF2 convex CH}.
It completes the proof.
\end{proof}


\subsection{Discrete energy dissipation law}

In our previous work \cite[Lemma 2.1]{LiaoZhang:2020linear},
the BDF2 kernels $b_{n-k}^{(n)}$ are shown to be positive definite
if the adjacent time-step ratios $0<r_k<\frac{3+\sqrt{17}}{2}$.
The following result shows that this sufficient condition
can be further improved in the theoretical manner. This improvement
is inspired by \cite[LemmaA.1]{LiaoJiZhang:2020pfc}
to find a lower bound for the eigenvalues of the step-scaled matrix $\widetilde{B}$,
see Lemma \ref{lem: tilde B-positiveDefinite} below.
For simplicity, we denote
\begin{align}\label{def: step-ratios function}
R_L\bra{z,s}:=\frac{2+4z-z^{3/2}}{1+z}-\frac{s^{3/2}}{1+s},
\quad\text{for $0< z,s< r_{*},$}
\end{align}
where $r_{*}\approx 4.864$ is the positive root of the equation $1+2r_{*}-r_{*}^{3/2}=0$.
According to the proof of \cite[LemmaA.1]{LiaoJiZhang:2020pfc},
$R_L\brab{z,s}$ is increasing in $(0,1)$
and decreasing in $(1, r_{*})$ with respect to $z$.
Also, $R_L\bra{z,s}$ is decreasing with respect to $s$ such that
\[
R_L\bra{z,s}>\min\{R_L\bra{0,r_{*}},R_L\brab{r_*,r_{*}}\}=\frac{2(1+2r_*-r_*^{3/2})}{1+z}=0
\quad\text{for $0< z,s<r_*$}.
\]

\begin{lemma}\label{lem: BDF2-Kernels-Positive}
Let $0<r_k<4.864$ for $2\le k\le N$. For any real sequence
$\{w_k\}_{k=1}^n$,
it holds that
\begin{align*}
2w_k\sum_{j=1}^kb_{k-j}^{(k)}w_j
\ge\frac{r_{k+1}^{3/2}}{1+r_{k+1}}\frac{w_k^2}{\tau_k}
-\frac{r_k^{3/2}}{1+r_k}\frac{w_{k-1}^2}{\tau_{k-1}}
+R_L(r_k,r_{k+1})
\frac{w_k^2}{\tau_k}\quad\text{for $k\ge2$}.
\end{align*}
So the discrete convolution kernels $b_{k-j}^{(k)}$ are positive definite in the sense that
\[
2\sum_{k=2}^n w_k \sum_{j=2}^k b_{k-j}^{(k)}w_j\ge \sum_{k=2}^nR_L(r_k,r_{k+1})\frac{w_k^2}{\tau_k}\quad\text{for $n\ge 2$}.
\]
\end{lemma}

\begin{proof}
Applying the inequality $-2ab\ge-a^2-b^2$, we take $u_k:=w_k/\sqrt{\tau_k}$ to find
\begin{align*}
2w_k\sum_{j=1}^k\frac{1}{\tau_k^2}b_{k-j}^{(k)}w_j
&=2\tau_kb_0^{(k)}u_k^2+2\sqrt{\tau_k\tau_{k-1}}b_1^{(k)}u_ku_{k-1}\\
&\ge \frac{2+4r_k}{1+r_k}u_k^2-\frac{r_k^{3/2}}{1+r_k}\bra{u_k^2+u_{k-1}^2}\\
&=\frac{r_{k+1}^{3/2}}{1+r_{k+1}}\frac{w_k^2}{\tau_k}
-\frac{r_k^{3/2}}{1+r_k}\frac{w_{k-1}^2}{\tau_{k-1}}
+R_L(r_k,r_{k+1})
\frac{w_k^2}{\tau_k}\quad\text{for $k\ge2$.}
\end{align*}
Summing this inequality from $k=2$ to $n$,
it is straightforward to obtain the claimed positive definiteness result.
It completes the proof.
\end{proof}

\begin{remark}
This lemma updates the sufficient condition of \cite[Lemma 2.1]{LiaoZhang:2020linear}.
Thus by following the discussions in \cite[Remark 3 and Remark 5]{LiaoZhang:2020linear},
one can verify that the variable-step BDF2 method is A-stable
if $0<r_k<4.864$ for $2\le k\le N$.
\end{remark}

Next theorem shows that
the numerical scheme \eqref{scheme: BDF2 convex CH} preserves a modified  energy dissipation property at the discrete levels,
and it is mesh-robustly stable in an energy norm.

\begin{theorem}\label{thm: energy-decay-law CH}
Let $\mathbf{S0}$ holds. If the stabilized parameter $A$
is properly large such that
\begin{align}\label{ieq: stabilized A restriction}
A\ge \frac{(r_{n}+r_{n+1}-1)^4}{64R_L^2(r_n,r_{n+1})}
\frac{\kappa^2}{\epsilon^2},
\end{align}
the convex-splitting BDF2 scheme \eqref{scheme: BDF2 convex CH}
preserves the following energy dissipation law
\begin{align*}
\mathcal{E}[\phi^n] \le \mathcal{E}[\phi^{n-1}]
\le \mathcal{E}[\phi^{1}]\quad\text{for $n\ge 2$,}
\end{align*}
where the modified discrete energy $\mathcal{E}[\phi^k]$  is defined by
\begin{align}\label{def: modified discrete energy}
\mathcal{E}[\phi^k]
:=E[\phi^k]
+\frac{\sqrt{\,r_{k+1}}\tau_{k+1}}{2\kappa(1+r_{k+1})}
\mynormb{\partial_{\tau} \phi^k}_{-1}^2
+\frac{\tau_k\tau_{k+1}}{2}\mynormb{\partial_{\tau}\phi^k}^2
+\frac{A\tau^2}{2}\mynormb{\nabla\phi^k}^2.
\end{align}
\end{theorem}
\begin{proof}
The volume conversation implies $\diff \phi^n\in\mathbb{\mathring V}_{h}$ for $n\ge1$.
Then we make the inner product of \eqref{scheme: BDF2 convex CH} by
$(-\Delta_h)^{-1}\diff \phi^{n}/\kappa$ and obtain
\begin{align}\label{Energy-Law-Inner}
\frac1\kappa\myinnerb{D_2\phi^n,\diff \phi^{n}}_{-1}
-\bra{\epsilon^2+A\tau^2}\myinnerb{\Delta_h\phi^n,\diff \phi^{n}}
+\myinnerb{\brat{\phi^n}^3-\hat{\phi}^n,\diff \phi^{n}}=0.
\end{align}
With the help of the summation by parts and $2a(a-b)=a^2-b^2+(a-b)^2$,
the second term at the left hand side of \eqref{Energy-Law-Inner} reads
\begin{align*}
\bra{\epsilon^2+A\tau^2}\myinnerb{\nabla_h\phi^n, \nabla_h\diff\phi^n}
=\frac{1}{2}\bra{\epsilon^2+A\tau^2}\brab{\mynormb{\nabla_h\phi^n}^2
-\mynormb{\nabla_h\phi^{n-1}}^2
+\mynormb{\nabla_h\diff \phi^{n}}^2}.
\end{align*}
It is easy to check the following identity
\begin{align*}
4a^{3}(a-b)=a^{4}-b^{4}+\kbra{2a^{2}+(a+b)^{2}}(a-b)^2.
\end{align*}
Then the nonlinear term in \eqref{Energy-Law-Inner} can be bounded by
\begin{align*}
\myinnerb{(\phi^{n})^3, \diff\phi^n}
\geq&\,\frac{1}{4}\mynormb{\phi^n}^{4}_{l^{4}}
-\frac{1}{4}\mynormb{\phi^{n-1}}^{4}_{l^{4}}.
\end{align*}
Noting the following identity
\begin{align*}
\kbra{(1+r_{n})b-r_{n}c}(a-b)
=&\,\frac{1}{2}(a^{2}-b^{2})+\frac{r_{n}}{2}(b-c)^{2}-\frac{r_{n+1}}{2}(a-b)^{2}\\
&\,+\frac{\bar{r}_n}{2}(a-b)^{2}-\frac{r_{n}}{2}(a-2b+c)^{2},
\end{align*}
where $\bar{r}_n:=r_{n}+r_{n+1}-1$ for brevity.
Then the extrapolation term in \eqref{Energy-Law-Inner} can be treated by
\begin{align*}
\myinnerb{\hat{\phi}^{n}, \diff\phi^n}=&\,\frac{1}{2}\brab{\mynormb{\phi^n}^{2}-\mynormb{\phi^{n-1}}^{2}}
+\frac{r_{n}}{2}\mynormb{\diff\phi^{n-1}}^{2}
-\frac{r_{n+1}}{2}\mynormb{\diff\phi^{n}}^{2}\nonumber\\
&\,+\frac{\bar{r}_n}{2}\mynormb{\diff\phi^{n}}^{2}
-\frac{r_{n}}{2}\mynormb{\diff\diff\phi^{n}}^{2}.
\end{align*}
The condition of \eqref{ieq: stabilized A restriction} gives that
$R_L(r_n,r_{n+1})\ge \kappa \bar{r}_n^2/\brat{8\epsilon A^{1/2}}.$
Taking $w_j=\diff \phi^{j}$ in the first inequality
of Lemma \ref{lem: BDF2-Kernels-Positive}, it is not difficult to get
\begin{align*}
\frac1\kappa\myinnerb{D_2\phi^n,\diff \phi^{n}}_{-1}
\ge&\,
\frac{\sqrt{\,r_{n+1}}\tau_{n+1}}{2\kappa(1+r_{n+1})}\mynormb{\partial_\tau\phi^n}_{-1}^2
-\frac{\sqrt{\,r_n}\tau_n}{2\kappa(1+r_n)}\mynormb{\partial_\tau\phi^{n-1}}_{-1}^2
+\frac{\bar{r}_n^2\epsilon^{-1}}{16\tau_nA^{\frac12}}\mynormb{\diff\phi^n}_{-1}^2.
\end{align*}

Thus it follows from \eqref{Energy-Law-Inner} that
\begin{align}\label{Energy-Inequality}
\mathcal{E}[\phi^n]
+\frac{\bar{r}_n^2\epsilon^{-1}}{16\tau_nA^{\frac12}}\mynormb{\diff\phi^n}_{-1}^2
+\frac{1}{2}\bra{\epsilon^2+A\tau^2}\mynormb{\nabla_h\diff \phi^{n}}^2
-\frac{\bar{r}_n}{2}\mynormb{\diff\phi^{n}}^{2}
\le \mathcal{E}[\phi^{n-1}]
\end{align}
for $n\ge 2$.
Recalling the definition of the maximum time-step $\tau$, one has
\[
\frac12\bra{\epsilon^2+A\tau^2}\mynormb{\nabla_h\diff \phi^{n}}^2
\ge\epsilon\tau A^{\frac12}\mynormb{\nabla_h\diff \phi^{n}}^2
\ge \epsilon\tau_nA^{\frac12}\mynormb{\nabla_h\diff \phi^{n}}^2.
\]
An application of the generalized H\"{o}lder inequality
\eqref{ieq: generalized Holder H-1} obtains
\begin{align*}
\frac{\bar{r}_n}{2}\mynormb{\diff \phi^n}^2
\le\frac{\abs{\bar{r}_n}}{2}\mynormb{\nabla_h\diff \phi^{n}}\mynormb{\diff \phi^n}_{-1}
\le \epsilon\tau_nA^{\frac12}\mynormb{\nabla_h\diff \phi^{n}}^2
+\frac{\bar{r}_n^2\epsilon^{-1}}{16\tau_nA^{\frac12}}\mynormb{\diff\phi^n}_{-1}^2.
\end{align*}
Combining it with \eqref{Energy-Inequality} yields
$\mathcal{E}[\phi^n] \le \mathcal{E}[\phi^{n-1}]$ for $n\ge 2$.
It completes the proof.
\end{proof}

\begin{remark}\label{remark: comments on stabilized A condition}
It is seen that this stabilization parameter constraint
\eqref{ieq: stabilized A restriction} requires $ A=O(\kappa^2/\epsilon^2)$.
Recalling the monotonicity of function $R_L(z,s)$,
we detail some requirements of $A$ to ensure energy stability:
\begin{itemize}[itemindent=0.1em]
\item[(i)] If time-step ratios $0<r_n,r_{n+1}\le 2$,
  and then $R_L(r_n,r_{n+1})\ge R_L(0,2)=2-\frac{2\sqrt{2}}{3}$.
  One needs  $ A\ge \frac{\bra{2+2-1}^4}{64R_L^2(0,2)}\frac{\kappa^2}{\epsilon^2}\approx 1.133\kappa^2/\epsilon^2$.
\item[(ii)] If time-step ratio $2<r_n\le 3$, one can choose $r_{n+1}$
  such that $0<r_{n+1}\le 2$,
   and then $R_L(r_n,r_{n+1})\ge R_L(3,2)=\frac{42-9\sqrt{3}-8\sqrt{2}}{12}$.
   It requires    $ A\ge 2.527\kappa^2/\epsilon^2$.
   \item[(iii)] If the current ratio is somewhat large such that $3<r_n\le r_{\text{user}}$, one can choose a small ratio $r_{n+1}$.
       For example, the step-ratio $r_{\text{user}}=4$ taken in adaptive time-steps computations shows that
        $ A\ge 1.778\kappa^2/\epsilon^2$
       is enough if $0<r_{n+1}\le 1$.
 \end{itemize}
 In Section 5, we consider the model parameters $\kappa=2\times10^{-3}$, $\epsilon=5\times10^{-2}$
 and $r_{\text{user}}=4$ for adaptive simulations.
 In such case, a mild constraint $ A\ge3/625$ is sufficient.
\end{remark}

\begin{remark}
The stabilized technique was originally introduced by Xu and Tang \cite{XuTang:2006Stability}
to build large time-stepping semi-implicit methods for phase filed models.
After that, various artificial stabilization terms
were proposed,
for instance, the second-order stabilization terms \cite{XuTang:2006Stability, WangYu:2018OnStabilizedCH}
$A\tau\Delta_h\bra{\phi^n-\phi^{n-1}}$, $A\tau\bra{\phi^n-\phi^{n-1}}$
and $A\bra{\phi^n-2\phi^{n-1}+\phi^{n-2}}$,
such that the discrete energy stability holds unconditionally
(or with reasonable stability condition);
however, the energy stability were all based on the assumptions that
nonlinear force $F'(\Phi)$ is Lipschitz continuous or the derivative of $F'(\Phi)$
is uniformly bounded.
He et al. \cite{HeLiuTang:2007StabilizedCH} used
the first-order stabilized term $A\Delta_h\bra{\phi^n-\phi^{n-1}}$
in which the energy stability
relayed on uniform bounds of the maximum norm of the numerical solutions.
Recently, under the time-step ratio $0<r_k<3+\sqrt{17}/2\approx3.561$,
the stabilized term $A\tau_n\Delta_h\bra{\phi^n-\phi^{n-1}}$
was first introduced in the variable-steps BDF2 method for CH model to achieve the unconditionally modified energy dissipation law \cite{ChenWangYanZhang:2019};
while the stabilization parameter $A$
could blow up for time-step ratios $r_k\rightarrow3.561.$
In current work,  under the time-step ratio condition $\mathbf{S0}$,
we introduce a new second-order stabilization term by adding a dissipation term $ A\tau^2\Delta_h\phi^n$ to ensure the energy stability for
the convex-splitting BDF2 scheme \eqref{scheme: BDF2 convex CH}.
Although the stabilized term $A\tau^2\Delta_h\phi^n$ is taken
as the maximum time step $\tau$ in every time step,
it avoids all the assumptions of nonlinear force $F'(\Phi)$
and the bounds of numerical solutions.
Meanwhile, the new artificial diffusion coefficient $A$
is bounded under the condition $\mathbf{S0}$.
Specially, the detailed discussion in Remark \ref{remark: comments on stabilized A condition}
shows that this new  artificial diffusion coefficient $A$
is of order $\kappa^2/\epsilon^2$ in the practical numerical computations.
\end{remark}

\begin{lemma}\label{lem: bound-BDF2 Solution CH}
Let $\mathbf{S0}$ and $\mathbf{A1}$ hold.
If the stabilized parameter $A$ fulfills \eqref{ieq: stabilized A restriction},
the solution of BDF2 time-stepping scheme \eqref{scheme: BDF2 convex CH}
is bounded in the sense that
\begin{align*}
\mynormb{\phi^n}+\mynormb{\nabla_h\phi^n}\le c_1:=\sqrt{4\epsilon^{-2}c_0+(2+\epsilon^2)\abs{\Omega_h}}\quad\text{for $n\ge2$,}
\end{align*}
where $c_1$ is dependent on the domain $\Omega$, the interface  parameter $\epsilon$ and the starting value $\phi^1$, but independent of the time $t_n$, the time-step sizes $\tau_n$ and the time-step ratios $r_n$.
\end{lemma}


\begin{proof}
Under the assumption $\mathbf{A1}$, the definition \eqref{def: modified discrete energy} of $\mathcal{E}[\phi^n]$ gives
\begin{align*}
\mathcal{E}[\phi^1]\le E[\phi^1]
+\frac{\tau_{2}}{2\kappa}\mynormb{\partial_{\tau} \phi^1}_{-1}^2
+\frac{\tau_1\tau_2}{2}\mynormb{\partial_{\tau} \phi^1}^2
+\frac{A\tau^2}{2}\mynormb{\nabla_h\phi^1}^2\le c_0.
\end{align*}
Thus the discrete energy dissipation law in Theorem \ref{thm: energy-decay-law CH} implies $c_0\ge \mathcal{E}[\phi^n]\ge E[\phi^{n}].$
Reminding the inequality $\mynorm{\phi^n}_{l^4}^4\ge 2(1+\epsilon^2)\mynorm{\phi^n}^2-(1+\epsilon^2)^2\abs{\Omega_h}$, due to the simple fact $(a^2-1-\epsilon^2)^2\ge0$,
one applies the definition \eqref{def: discrete free energy} of $E[\phi^{n}]$ to get
\begin{align*}
4c_0\ge&\, 2\epsilon^2\mynormb{\nabla_h\phi^n}^2+4\myinnerb{F(\phi^n),1}
\ge 2\epsilon^2\mynormb{\nabla_h\phi^n}^2+2\epsilon^2\mynorm{\phi^n}^2-\epsilon^2(2+\epsilon^2)\abs{\Omega_h},
\end{align*}
and then
\begin{align*}
\brab{\mynormb{\phi^n}+\mynormb{\nabla_h\phi^n}}^2\le
2\mynorm{\phi^n}^2+2\mynormb{\nabla_h\phi^n}^2\le 4\epsilon^{-2}c_0+(2+\epsilon^2)\abs{\Omega_h}\quad\text{for $n\ge2$.}
\end{align*}
It implies the claimed result and completes the proof.
\end{proof}

\section{Some discrete convolution inequalities}
\setcounter{equation}{0}

Our error analysis is closely related to the convolution form
\eqref{scheme: DOC form BDF2 convex CH}, so we
need some detail properties and discrete convolution inequalities with respect to the DOC kernels $\theta^{(n)}_{n-j}$.
It is to emphasize that the positive constants $\mathfrak{m}_1$, $\mathfrak{m}_2$ and $\mathfrak{m}_3$ involved in this section
are independent of the time $t_n$, time-step sizes $\tau_n$ and the step ratios $r_n$.
Actually, they would take different values for different choices of step ratios $r_n$,
but are bounded with respect to the changes of step ratios,
even when $r_n$ approaches the user limit $r_{\mathrm{user}}$.

\subsection{Simple properties of DOC kernels}
Following the proofs of \cite[Lemma 2.2, Corollary 2.1 and Lemma 2.3]{LiaoZhang:2020linear},
we can obtain some simple properties of  the DOC kernels.

\begin{lemma}\label{lem: DOC property}
	If $\mathbf{S0}$ holds, the DOC kernels $\theta_{n-j}^{(n)}$ defined in \eqref{def: DOC-Kernels} satisfy:
\begin{itemize}
  \item[(I)] The discrete kernels $\theta_{n-j}^{(n)}$ are positive definite;
  \item[(II)] The discrete kernels $\theta_{n-j}^{(n)}$ are positive and
  $\displaystyle \theta_{n-j}^{(n)}=\frac{1}{b^{(j)}_{0}}\prod_{i=j+1}^n\frac{r_i^2}{1+2r_i}$
  for $2\le j\le n$;
  \item[(III)] $\displaystyle \sum_{j=2}^{n}\theta_{n-j}^{(n)}\le\tau_n$ such that
  $\displaystyle \sum_{k=2}^{n}\sum_{j=2}^{k}\theta_{k-j}^{(k)}\le t_n$ for $n\ge2$.
\end{itemize}
\end{lemma}

We introduce the following two $(n-1)\times (n-1)$ matrices
\[
 B_2:=
\left(
\begin{array}{cccc}
b_{0}^{(2)}  &             &            & \\
b_{1}^{(3)} &b_{0}^{(3)}  &            & \\
            &\ddots       &\ddots      &\\
            &             &b_{1}^{(n)} &b_{0}^{(n)}  \\
\end{array}
\right)\quad\text{and}\quad
\Theta_2:=
\left(
\begin{array}{cccc}
\theta_{0}^{(2)}  &                  &  & \\
\theta_{1}^{(3)}  &\theta_{0}^{(3)}  &  & \\
\vdots           &\vdots           &\ddots  &\\
\theta_{n-2}^{(n)}&\theta_{n-3}^{(n)}&\cdots  &\theta_{0}^{(n)}  \\
\end{array}
\right),
\]
where the discrete kernels
$b_{n-k}^{(n)}$ and $\theta_{n-k}^{(n)}$ are defined by \eqref{def: BDF2-kernels} and \eqref{def: DOC-Kernels}, respectively.
It follows from the discrete orthogonal
identity  \eqref{eq: orthogonal identity} that
\begin{align}\label{matrix: orthogonal identity}
\Theta_2= B_2^{-1}.
\end{align}
If the step ratios condition $\mathbf{S0}$ holds, Lemma \ref{lem: BDF2-Kernels-Positive} shows that
the real symmetric matrix
\begin{align}\label{matrix: B}
 B:= B_2+ B_2^T
\end{align}
is positive definite, that is,
\begin{align*}
\myvec{w}^T B\myvec{w}=2\sum_{k=2}^n w^k \sum_{j=2}^k b_{k-j}^{(k)}w^j\ge
\sum_{k=2}^n
\frac{R_L(r_k,r_{k+1})}{\tau_k}(w^k)^2\,,
\end{align*}
where the function $R_L(z,s)$ is defined by \eqref{def: step-ratios function} and the vector $\myvec{w}:=(w^2,w^3,\cdots,w^n)^T$.
According to Lemma \ref{lem: DOC property} (I), the following symmetric matrix
\begin{align}\label{matrix: Theta}
\Theta:=\Theta_2+\Theta_2^T
= B_2^{-1}+( B_2^{-1})^T
=( B_2^{-1})^T(B_2+B_2^T)B_2^{-1}=(B_2^{-1})^TBB_2^{-1}
\end{align}
is also positive definite in the sense of $\myvec{w}^T\Theta\myvec{w}=2\sum_{k,j}^{n,k}\theta_{k-j}^{(k)}w^jw^k>0$.
Here and hereafter, we denote
$\sum_{k,j}^{n,k}:=\sum_{k=2}^n\sum_{j=2}^k$ for the simplicity of presentation.

\subsection{Eigenvalue estimates}

To facilitate the proofs in what follows, we are to define the following step-scaled matrix
\begin{align}\label{matrix: modified-B2}
\widetilde{B}_2:=\Lambda_\tau B_2\Lambda_\tau=
\left(
\begin{array}{cccc}
\tilde{b}_{0}^{(2)}  &             &            & \\
\tilde{b}_{1}^{(3)} &\tilde{b}_{0}^{(3)}  &            & \\
            &\ddots       &\ddots      &\\
            &             &\tilde{b}_{1}^{(n)} &\tilde{b}_{0}^{(n)}  \\
\end{array}
\right)_{(n-1)\times(n-1)},
\end{align}
where the diagonal matrix
$\Lambda_\tau:=\text{diag}\bra{\sqrt{\tau_2},\sqrt{\tau_3},\cdots,\sqrt{\tau_n}}$
so that the step-scaled discrete kernels $\tilde{b}_{0}^{(k)}$ and $\tilde{b}_{1}^{(k)}$ are given by
\begin{align}\label{def: tilde kernels b}
\tilde{b}_{0}^{(k)}=\frac{1+2r_k}{1+r_k}\quad
\text{and}\quad
\tilde{b}_{1}^{(k)}=-\frac{r_k^{3/2}}{1+r_k}
\quad\text{for $2\leq k\leq n$}.
\end{align}
Moreover, we will use the following real symmetric  matrix,
\begin{align}\label{matrix: two modified-B}
\widetilde{B}:=\widetilde{B}_2+\widetilde{B}_2^T=\Lambda_\tau B\Lambda_\tau.
\end{align}

The following two lemmas present some eigenvalue estimates of $\widetilde{B}$ and $\widetilde{B}_2^T\widetilde{B}_2$.
To avoid possible confusions,
we define the vector norm $\timenorm{\cdot}$ by
$\timenorm{\myvec{u}}:=\sqrt{\myvec{u}^T\myvec{u}}$ for any real vector $\myvec{u}$ and
the associated matrix norm
$\timenorm{U}:=\sqrt{\lambda_{\max}\brab{U^TU}}$.

\begin{lemma}\label{lem: tilde B-positiveDefinite}
If $\mathbf{S0}$ holds, there exists a positive constant $\mathfrak{m}_1$
 such that $\lambda_{\min}\brab{\widetilde{B}}\ge \mathfrak{m}_1>0$.
\end{lemma}

\begin{proof}
This proof can be followed from \cite[Lemma A.1]{LiaoJiZhang:2020pfc}. We include the main ingredient for the completeness.
Applying the Gerschgorin's circle theorem to the matrix $\widetilde{B}$, one has
\begin{align*}
\lambda_{\min}\brab{\widetilde{B}}
\ge\min_{2\leq k\leq n}R_L\bra{r_k,r_{k+1}}
>R_L\bra{r_{\mathrm{user}},r_{\mathrm{user}}}
=\frac{2(1+2r_{\mathrm{user}}-r_{\mathrm{user}}^{3/2})}{1+r_{\mathrm{user}}}>0,
\end{align*}
where $R_L\brab{z,s}$ is defined by \eqref{def: step-ratios function}.
It completes the proof by taking $\mathfrak{m}_1=\frac{2(1+2r_{\mathrm{user}}-r_{\mathrm{user}}^{3/2})}{1+r_{\mathrm{user}}}.$
\end{proof}

\begin{lemma}\label{lem: tilde B2TB2-positiveDefinite}
If $\mathbf{S0}$ holds, there exists a positive constant  $\mathfrak{m}_2$
 such that $\lambda_{\max}\brab{\widetilde{B}_2^T\widetilde{B}_2}\leq \mathfrak{m}_2$.
\end{lemma}
\begin{proof}
This proof can be followed from \cite[Lemma A.2]{LiaoJiZhang:2020pfc}. We include the main ingredient for the completeness.
By writing out the tri-diagonal matrix $\widetilde{B}_2^T\widetilde{B}_2$ and applying the Gerschgorin's circle theorem, one can find
\begin{align*}
\lambda_{\max}\brab{\widetilde{B}_2^T\widetilde{B}_2}\le\max_{2\le k\le n}R_U\bra{r_{k},r_{k+1}}<R_U\bra{r_{\mathrm{user}},r_{\mathrm{user}}},
\end{align*}
where the function $R_U\bra{z,s}$ is defined by
\begin{align*}
R_U\bra{z,s}:=\frac{(1+2z)(1+2z+z^{3/2})}{\bra{1+z}^2}
+\frac{s^{3/2}(1+2s+s^{3/2})}{\bra{1+s}^2}
\quad\text{for $0\le z,s< r_{\mathrm{user}}.$}
\end{align*}
An upper bound is then obtained by taking $\mathfrak{m}_2=R_U\bra{r_{\mathrm{user}},r_{\mathrm{user}}}$.
\end{proof}


By the above two lemmas, we can bound the minimum eigenvalue of $\Theta$.
\begin{lemma}\label{lem: Theta minimum eigenvalue}
If $\mathbf{S0}$ holds, the real symmetric  matrix $\Theta$ in \eqref{matrix: Theta} satisfies
\begin{align*}
\myvec{v}^T\Theta\myvec{v}
\geq&\,\frac{\mathfrak{m}_1}{\mathfrak{m}_2}\,\timenorm{\Lambda_\tau\myvec{v}}^2\quad \text{for any vector $\myvec{v}$}.
\end{align*}
\end{lemma}
\begin{proof}
Lemma \ref{lem: tilde B-positiveDefinite} says that real
symmetric  matrix $\widetilde{B}$ is positive definite.
There exists a non-singular upper triangular matrix $\widetilde{U}$ such that $\widetilde{B}=\widetilde{U}^T\widetilde{U}$.
By using \eqref{matrix: Theta} and \eqref{matrix: two modified-B}, one gets
\begin{align*}
\myvec{v}^T\Theta\myvec{v}=\myvec{v}^T(B_2^{-1})^TBB_2^{-1}\myvec{v}
=\myvec{v}^T(B_2^{-1})^T\Lambda_\tau^{-1}\widetilde{B}\Lambda_\tau^{-1}B_2^{-1}\myvec{v}
=\timenorm{\widetilde{U}\Lambda_\tau^{-1}B_2^{-1}\myvec{v}}^2.
\end{align*}
Thus it follows that
\begin{align*}
\timenorm{\Lambda_\tau\myvec{v}}^2=&\,\timenorm{\Lambda_\tau B_2\Lambda_\tau\widetilde{U}^{-1}\widetilde{U}\Lambda_\tau^{-1}B_2^{-1}\myvec{v}}^2
\leq\timenorm{\widetilde{B}_2\widetilde{U}^{-1}}^2\timenorm{\widetilde{U}\Lambda_\tau^{-1}B_2^{-1}\myvec{v}}^2\\
\leq&\,\timenorm{\widetilde{B}_2}^2\timenorm{\widetilde{U}^{-1}}^2\myvec{v}^T\Theta\myvec{v}
=\lambda_{\max}\brab{\widetilde{B}_2^T\widetilde{B}_2}\lambda_{\max}\brab{\widetilde{B}^{-1}}\myvec{v}^T\Theta\myvec{v}.
\end{align*}
Thus Lemmas \ref{lem: tilde B-positiveDefinite}
and \ref{lem: tilde B2TB2-positiveDefinite} yield the claimed inequality.
\end{proof}

To evaluate the maximum eigenvalue of $\Theta$,
consider the inverse matrix of the matrix $\widetilde{B}_2$,
\begin{align}\label{matrix: modified-Theta2}
\widetilde{\Theta}_2:=\widetilde{B}_2^{-1}=\Lambda_\tau^{-1}\Theta_2\Lambda_\tau^{-1}=
\left(
\begin{array}{cccc}
\tilde{\theta}_{0}^{(2)}  &                  &  & \\
\tilde{\theta}_{1}^{(3)}  &\tilde{\theta}_{0}^{(3)}  &  & \\
\vdots           &\vdots           &\ddots  &\\
\tilde{\theta}_{n-2}^{(n)}&\tilde{\theta}_{n-3}^{(n)}&\cdots  &\tilde{\theta}_{0}^{(n)}  \\
\end{array}
\right),
\end{align}
where the step-scaled DOC kernels $\tilde{\theta}_{k-j}^{(k)}$
follow from Lemma \ref{lem: DOC property} (II),
\begin{align}\label{def: tilde kernels theta}
\tilde{\theta}_{k-j}^{(k)}:=\frac1{\sqrt{\tau_k\tau_j}}\theta_{k-j}^{(k)}=\frac{1+r_j}{1+2r_j}\prod_{i=j+1}^k\frac{r_i^{3/2}}{1+2r_i}
\quad\text{for $2\leq j\le k\leq n$}.
\end{align}

\begin{lemma}\label{lem: Theta maximum eigenvalue}
If $\mathbf{S0}$ holds, then there exists a positive constant  $\mathfrak{m}_3$
 such that
\begin{align*}
\myvec{v}^T\Theta\myvec{v}
\leq&\,\mathfrak{m}_3\,\timenorm{\Lambda_\tau\myvec{v}}^2\quad \text{for any vector $\myvec{v}$}.
\end{align*}
\end{lemma}
\begin{proof}
Let $\widetilde{\Theta}=\widetilde{\Theta}_2+\widetilde{\Theta}_2^T$.
Since $0<\frac{x^{3/2}}{1+2x}<m_{*}:=\frac{r_{\mathrm{user}}^{3/2}}{1+2r_{\mathrm{user}}}<1$ for any $x\in [0,r_{\mathrm{user}}]$,
one can apply the formula \eqref{def: tilde kernels theta} to get
\begin{align*}
\mathfrak{R}_{n,k}:=\sum_{j=2}^k\tilde{\theta}_{k-j}^{(k)}+\sum_{j=k}^n\tilde{\theta}_{j-k}^{(j)}
\le \sum_{j=2}^km_{*}^{k-j}+\sum_{j=k}^nm_{*}^{j-k}<\frac{2}{1-m_{*}}\quad\text{for $2\le k\le n$.}
\end{align*}
One has $\lambda_{\max}\brab{\widetilde{\Theta}}\le\max_{2\le k\le n}\mathfrak{R}_{n,k}<\mathfrak{m}_3:=\frac{2}{1-m_{*}}$ by the Gerschgorin's circle
theorem.
It implies $\myvec{w}^T\widetilde{\Theta}\myvec{w}
\leq\mathfrak{m}_3\timenorm{\myvec{w}}^2$ for any $\myvec{w}$ and
the choice $\myvec{w}:=\Lambda_\tau\myvec{v}$ completes the proof.
\end{proof}

\subsection{Discrete convolution inequalities}

The following two lemmas describe the Young-type convolution inequality.

\begin{lemma}\label{lem:DOC quadr form Young inequ2}
If $\mathbf{S0}$ holds, then for any real sequences $\{v^k\}_{k=2}^n$ and $\{w^k\}_{k=2}^n$,
\begin{align*}
\sum_{k,j}^{n,k}\theta_{k-j}^{(k)} w^k v^j
\le\varepsilon\sum_{k,j}^{n,k}\theta_{k-j}^{(k)} v^k v^j
+\frac{1}{2\mathfrak{m}_1\varepsilon}
\sum_{k=2}^{n}\tau_k (w^k)^2
\quad \text{for $\forall\;\varepsilon > 0$.}
\end{align*}
\end{lemma}
\begin{proof}
Let $\myvec{w}:=\brat{w^2,w^3,\cdots,w^n}^T$. A similar proof of \cite[Lemma A.3]{LiaoJiZhang:2020pfc} gives
\begin{align*}
\sum_{k,j}^{n,k}\theta_{k-j}^{(k)}v^j w^k
\le \varepsilon\sum_{k,j}^{n,k}\theta_{k-j}^{(k)}v^j v^k
+\frac{1}{2\varepsilon}\myvec{w}^T B^{-1}\myvec{w}
\quad \text{for any $\varepsilon > 0$}.
\end{align*}
From the proof Lemma \ref{lem: Theta minimum eigenvalue}, we have
$B^{-1}=\Lambda_\tau \widetilde{U}^{-1}\brab{\Lambda_\tau \widetilde{U}^{-1}}^T$ and then
\begin{align*}
\myvec{w}^T B^{-1}\myvec{w}=&\,\myvec{w}^T\Lambda_\tau \widetilde{U}^{-1}
\brab{\Lambda_\tau \widetilde{U}^{-1}}^T\myvec{w}
=\timenorm{\brab{\widetilde{U}^{-1}}^T\Lambda_\tau \myvec{w}}^2\\
\leq&\,\timenorm{\brab{\widetilde{U}^{-1}}^T}^2
\timenorm{\Lambda_\tau\myvec{w}}^2
=\lambda_{\max}\brab{(\widetilde{B})^{-1}}
\myvec{w}^T\Lambda_\tau^2\myvec{w}\leq \mathfrak{m}_1^{-1}\sum_{k=2}^{n}\tau_k (w^k)^2,
\end{align*}
where Lemma \ref{lem: tilde B-positiveDefinite} has been used.
It completes the proof.
\end{proof}

\begin{lemma}\label{lem:DOC quadr form Young inequ-embedding}
If $\mathbf{S0}$ holds, then for any real sequences $\{v^k\}_{k=2}^n$ and $\{w^k\}_{k=2}^n$,
\begin{align*}
\sum_{k,j}^{n,k}\theta_{k-j}^{(k)} w^k v^j
\le\varepsilon\sum_{k=2}^n\tau_k  (v^k)^2
+\frac{\mathfrak{m}_3}{4\mathfrak{m}_1\varepsilon}
\sum_{k=2}^{n}\tau_k (w^k)^2\quad \text{for $\forall\;\varepsilon > 0$}.
\end{align*}
\end{lemma}

\begin{proof}For fixed time index $n$, taking $\varepsilon:=2\varepsilon_0/\mathfrak{m}_3$ in Lemma \ref{lem:DOC quadr form Young inequ2} yields
\begin{align*}
\sum_{k,j}^{n,k}\theta_{k-j}^{(k)} w^k v^j
\le&\,\frac{2\varepsilon_0}{\mathfrak{m}_3}\sum_{k,j}^{n,k}\theta_{k-j}^{(k)} v^k v^j
+\frac{\mathfrak{m}_3}{4\mathfrak{m}_1\varepsilon_0}\sum_{k=2}^{n}\tau_k (w^k)^2\\
\le&\,\varepsilon_0\sum_{k=2}^{n}\tau_k (v^k)^2
+\frac{\mathfrak{m}_3}{4\mathfrak{m}_1\varepsilon_0}
\sum_{k=2}^{n}\tau_k (w^k)^2,
\end{align*}
where Lemma \ref{lem: Theta maximum eigenvalue} was used in the last inequality.
It completes the proof by choosing $\varepsilon_0:=\varepsilon$.
\end{proof}

We now present two discrete  embedding-type convolution inequalities
by considering three time-space discrete functions $u^k$, $v^k$ and $w^k$ $(2\le k\le n)$
in the space $\mathbb{V}_{h}$ or its subspace $\mathbb{\mathring V}_{h}$.

\begin{lemma}\label{lem:DOC quadr form H1 embedding inequ}
Assume that $u^k, w^k\in \mathbb{V}_{h}$, $v^k\in \mathbb{\mathring V}_{h}$ $(2\le k\le n)$
and there exists a constant $c_u$ such that
$\mynormb{u^k}_{l^3}\le c_u$ for $2\le k\le n$. If $\mathbf{S0}$ holds, then for any $\varepsilon >0$,
\begin{align*}
\sum_{k,j}^{n,k}\theta_{k-j}^{(k)} \myinnerb{u^jv^j,w^k}
\le
\varepsilon\sum_{k,j}^{n,k}\theta_{k-j}^{(k)}\myinnerb{\nabla_h v^j,\nabla_h v^k}
+\frac{c_z^2c_u^2\mathfrak{m}_2\mathfrak{m}_3}{2\mathfrak{m}_1^2\varepsilon}
\sum_{k=2}^{n}\tau_k \mynormb{w^k}^2.
\end{align*}
\end{lemma}
\begin{proof}For fixed time index $n$, taking $v^j:=u_h^jv_h^j$
and $\varepsilon:=\varepsilon_1$ in Lemma \ref{lem:DOC quadr form Young inequ-embedding}, we have
\begin{align*}
\sum_{k,j}^{n,k}\theta_{k-j}^{(k)} \myinnerb{u^jv^j,w^k}
\le
\varepsilon_1\sum_{k=2}^n\tau_k  \mynormb{u^kv^k}^2
+\frac{\mathfrak{m}_3}{4\mathfrak{m}_1\varepsilon_1}\sum_{k=2}^{n}\tau_k \mynormb{w^k}^2.
\end{align*}
The well--known H\"{o}lder inequality and the discrete embedding inequality \eqref{ieq: mean-zero H1 embedding L6} imply that
$\mynormb{u^kv^k}\le \mynormb{u^k}_{l^3}\mynormb{v^k}_{l^6}\le c_z\mynormb{u^k}_{l^3}\mynormb{\nabla_h v^k}
\le c_zc_u\mynormb{\nabla_h v^k}$. We derive that
\begin{align*}
\sum_{k=2}^n\tau_k  \mynormb{u^kv^k}^2\le&\, c_z^2c_u^2\sum_{k=2}^n\tau_k\mynormb{\nabla_h v^k}^2.
\end{align*}
Then it follows that
\begin{align}\label{lemProof:DOC quadr form H1 embedding inequ}
\sum_{k,j}^{n,k}\theta_{k-j}^{(k)} \myinnerb{u^jv^j,w^k}
\le\varepsilon_1c_z^2c_u^2\sum_{k=2}^n\tau_k\mynormb{\nabla_h v^k}^2
+\frac{\mathfrak{m}_3}{4\mathfrak{m}_1\varepsilon_1}\sum_{k=2}^{n}\tau_k \mynormb{w^k}^2.
\end{align}
Following the proof of Lemma \ref{lem: Theta minimum eigenvalue},
it is not difficult to get (cf. \cite{LiaoSongTangZhou:2020mbe})
\begin{align*}
\sum_{k=2}^n\tau_k\mynormb{\nabla_h v^k}^2\le\frac{2\mathfrak{m}_2}{\mathfrak{m}_1}
\sum_{k,j}^{n,k}\theta_{k-j}^{(k)}\myinnerb{\nabla_h v^j,\nabla_h v^k}.
\end{align*}
Inserting this inequality into \eqref{lemProof:DOC quadr form H1 embedding inequ} and choosing the parameter
$\varepsilon_1:=\mathfrak{m}_1\varepsilon/\brat{2c_z^2c_u^2\mathfrak{m}_2},$
we get the claimed inequality and complete the proof.
\end{proof}

\begin{lemma}\label{lem:DOC quadr form H2-H2 embedding inequ}
Assume that $u^k\in \mathbb{V}_{h}$, $w^k\in \mathbb{\mathring V}_{h}$ $(2\le k\le n)$
and there exists a constant $c_u$ such that
$\mynormb{u^k}_{l^3}\le c_u$ for $2\le k\le n$. If $\mathbf{S0}$ holds, then for any $\varepsilon >0$,
\begin{align*}
\sum_{k,j}^{n,k}\theta_{k-j}^{(k)} \myinnerb{u^jw^j,\Delta_h w^k}
\le\varepsilon\sum_{k,j}^{n,k}\theta_{k-j}^{(k)}\myinnerb{\Delta_h w^j,\Delta_h w^k}
+\frac{c_z^4c_u^4\mathfrak{m}_2^3\mathfrak{m}_3^2}{\mathfrak{m}_1^5\varepsilon^3}
\sum_{k=2}^n\tau_k  \mynormb{w^k}^2.
\end{align*}
\end{lemma}

\begin{proof}For fixed time index $n$,   we start the proof from \eqref{lemProof:DOC quadr form H1 embedding inequ}
by setting $w^j:=\Delta_h w^j$, $v^j:=w^j$ and $\varepsilon_1:=\mathfrak{m}_2\mathfrak{m}_3/(\varepsilon_4\mathfrak{m}_1^2)$, that is,
\begin{align}\label{lemProof:DOC quadr form H2-H2 embedding inequ}
\sum_{k,j}^{n,k}\theta_{k-j}^{(k)} \myinnerb{u^jw^j,\Delta_h w^k}
\le&\,\frac{c_z^2c_u^2\mathfrak{m}_2\mathfrak{m}_3}{2\mathfrak{m}_1^2\varepsilon_4}
\sum_{k=2}^n\tau_k\mynormb{\nabla_h w^k}^2
+\frac{\mathfrak{m}_1\varepsilon_4}{2\mathfrak{m}_2}\sum_{k=2}^{n}\tau_k \mynormb{\Delta_h w^k}^2\nonumber\\
\le&\,\frac{c_z^2c_u^2\mathfrak{m}_2\mathfrak{m}_3}{2\mathfrak{m}_1^2\varepsilon_4}
\sum_{k=2}^n\tau_k\mynormb{\nabla_h w^k}^2
+\varepsilon_4\sum_{k,j}^{n,k}\theta_{k-j}^{(k)}\myinnerb{\Delta_h w^j,\Delta_h w^k},
\end{align}
where Lemma \ref{lem: Theta minimum eigenvalue} has been used to handle the last term.
Furthermore, by using the classical Young's inequality and Lemma \ref{lem: Theta minimum eigenvalue}, one gets
\begin{align*}
\sum_{k=2}^n\tau_k\mynormb{\nabla_h w^k}^2
=&\,\sum_{k=2}^n\tau_k \myinnerb{-\Delta_h w^k,w^k}
\le\frac{\varepsilon_3}{2}\sum_{k=2}^n\tau_k \mynormb{\Delta_h w^k}^2
+\frac{1}{2\varepsilon_3}\sum_{k=2}^n\tau_k  \mynormb{w^k}^2\nonumber\\
\le&\, \frac{\mathfrak{m}_2\varepsilon_3}{\mathfrak{m}_1}
\sum_{k,j}^{n,k}\theta_{k-j}^{(k)}\myinnerb{\Delta_h w^j,\Delta_h w^k}
+\frac{1}{2\varepsilon_3}\sum_{k=2}^n\tau_k  \mynormb{w^k}^2.
\end{align*}
Inserting this inequality into \eqref{lemProof:DOC quadr form H2-H2 embedding inequ}, we have
\begin{align*}
\sum_{k,j}^{n,k}\theta_{k-j}^{(k)} \myinnerb{u^jw^j,\Delta_h w^k}
\le&\,\braB{\frac{c_z^2c_u^2\mathfrak{m}_2^2\mathfrak{m}_3\varepsilon_3}{2\mathfrak{m}_1^3\varepsilon_4}+\varepsilon_4}
\sum_{k,j}^{n,k}\theta_{k-j}^{(k)}\myinnerb{\Delta_h w^j,\Delta_h w^k}\\
&\,+\frac{c_z^2c_u^2\mathfrak{m}_2\mathfrak{m}_3}{4\mathfrak{m}_1^2\varepsilon_3\varepsilon_4}
\sum_{k=2}^n\tau_k  \mynormb{w^k}^2.
\end{align*}
Now by choosing $\varepsilon_4:=\varepsilon/2$ and $\varepsilon_3:=\mathfrak{m}_1^3\varepsilon_4\varepsilon/(c_z^2c_u^2\mathfrak{m}_2^2\mathfrak{m}_3)$,
we obtain the claimed inequality.
\end{proof}

\section{Robust $L^2$ norm error estimate}
\setcounter{equation}{0}

\subsection{Convolutional consistency and technical lemma}


Let $\xi_{\Phi}^j$ be the local consistency errors of
the convex-splitting BDF2 scheme \eqref{scheme: BDF2 convex CH},
arising from the BDF2 formula \eqref{def: BDF2-Formula},
the extrapolation approximation and the artificial stabilization term,
at the time $t=t_j$, that is,
\begin{align}\label{def: BDF2-local consistency}
\xi_{\Phi}^j:=\kbrab{D_2\Phi(t_j)-\partial_t\Phi(t_j)}
+\kappa\kbrab{\Delta\hat{\Phi}(t_j)-\Delta\Phi(t_j)}
+\kappa A\tau^2\Delta^2\Phi(t_j).
\end{align}
We will consider a convolutional consistency error $\Xi_{\Phi}^k$ defined by
\begin{align}\label{def: BDF2-global consistency}
\Xi_{\Phi}^k:=\sum_{j=2}^k\theta_{k-j}^{(k)}\xi_{\Phi}^j
\quad\text{for $k\ge2$.}
\end{align}

\begin{lemma}\label{lem: BDF2-Consistency-Error}
If $\mathbf{S0}$ holds, the convolutional consistency error $\Xi_{\Phi}^k$ in \eqref{def: BDF2-global consistency} satisfies
\begin{align*}
\sum_{k=2}^n\absb{\Xi_{\Phi}^k}
\le&\, t_n\tau^2\max_{1\le j\le n}
\braB{3\absb{\Phi'''(t_j)}+2\kappa\absb{\Delta\Phi''(t_j)}
+\kappa A\absb{\Delta^2\Phi(t_j)}}\quad\text{for $n\ge2$.}
\end{align*}
\end{lemma}
\begin{proof}
By following the proof of \cite[Lemma 3.4]{LiaoJiZhang:2020pfc},
the convolution consistency error
for the BDF2 formula \eqref{def: BDF2-Formula} can be bounded by
\begin{align*}
\sum_{j=2}^k\theta_{k-j}^{(k)}
\absb{D_2\Phi(t_j)-\partial_t\Phi(t_j)}
\le 3\sum_{j=1}^{k}\theta_{k-j}^{(k)}
\tau_{j}\int_{t_{j-1}}^{t_j} \absb{\Phi'''(s)}\zd s \zd{s}\quad\text{for $k\ge2$.}
\end{align*}
By using the Taylor's expansion formula, one has
\begin{align*}
\hat{v}^j-v^j&=\int_{t_{j-1}}^{t_j}(s-t_j)v''(s)\zd s
-r_j\int_{t_{j-2}}^{t_{j-1}}(s-t_{j-2})v''(s)\zd s,
\end{align*}
which in turn yields (by taking $v:=\Delta\Phi$)
\begin{align*}
\sum_{j=2}^k\theta_{k-j}^{(k)}
\kappa\absb{\Delta\hat{\Phi}(t_j)-\Delta\Phi(t_j)}
\le \kappa\sum_{j=1}^{k}\theta_{k-j}^{(k)}
\tau_{j}\int_{t_{j-2}}^{t_j} \absb{\Delta\Phi''(s)}\zd s
\quad\text{for $k\ge2$.}
\end{align*}
For the stabilization term, it is straightforward to derive
\begin{align*}
\sum_{j=2}^k\theta_{k-j}^{(k)}
\absb{\kappa A\tau^2\Delta^2\Phi(t_j)}
\le \kappa A\tau^2\sum_{j=1}^{k}\theta_{k-j}^{(k)}
\absb{\Delta^2\Phi(t_j)}\quad\text{for $k\ge2$.}
\end{align*}
Collecting the above estimates and using Lemma \ref{lem: DOC property} (III),
one obtains the claimed estimate on the convolutional consistency  immediately.
This completes the proof.
\end{proof}

We use  the standard seminorms and norms in the Sobolev space $H^{m}(\Omega)$ for $m\ge0$.
Let $C_{per}^{\infty}(\Omega)$ be a set of infinitely differentiable $L$-periodic functions defined on $\Omega$,
and $H_{per}^{m}(\Omega)$ be the closure of $C_{per}^{\infty}(\Omega)$ in $H^{m}(\Omega)$,
endowed with the semi-norm $|\cdot|_{H_{per}^m}$ and the norm $\mynorm{\cdot}_{H_{per}^{m}}$.

For simplicity, denote $|\cdot|_{H^m}:=|\cdot|_{H_{per}^m}$, $\mynorm{\cdot}_{H^{m}}:=\mynorm{\cdot}_{H_{per}^{m}}$, and $\mynorm{\cdot}_{L^{2}}:=\mynorm{\cdot}_{H^{0}}$.
Next lemma lists some approximations, cf. \cite{ShenTangWang:2011Spectral,ShenWangWangWise:2012},
of the $L^2$-projection operator $P_{M}$ and
trigonometric interpolation operator $I_{M}$ defined in subsection 2.1.
\begin{lemma}\label{lem:Projection-Estimate}
For any $u\in{H_{per}^{q}}(\Omega)$ and $0\le{s}\le{q}$, it holds that
\begin{align}
\mynorm{P_{M}u-u}_{H^{s}}
\le C_uh^{q-s}|u|_{H^{q}},
\quad \mynorm{P_{M}u}_{H^{s}}\le C_u\mynorm{u}_{H^{s}};\label{Projection-Estimate}
\end{align}
and, in addition if $q>3/2$,
\begin{align}
\mynorm{I_{M}u-u}_{H^{s}}
\le C_uh^{q-s}|u|_{H^{q}},
\quad \mynorm{I_{M}u}_{H^{s}}\le C_u\mynorm{u}_{H^{s}}.\label{Interpolation-Estimate}
\end{align}
\end{lemma}

\subsection{Convergence analysis}

Note that,  the energy dissipation law \eqref{cont:energy dissipation}
of CH model \eqref{cont: Problem-CH} shows that $E[\Phi^n]\le E[\Phi(t_0)]$.
From the formulation \eqref{cont:free energy}, it is easy to check that $\mynormb{\Phi^n}_{H^1}$ can be bounded
by a time-independent constant.
Let $\Phi_M^n:=\brab{P_M\Phi}(\cdot,t_n)$
be the $L^2$-projection of exact solution at time $t=t_n$.
The projection estimate \eqref{Projection-Estimate} in Lemma \ref{lem:Projection-Estimate} yields
\begin{align}\label{H1 bound projection}
\mynormb{\Phi_M^n}+\mynormb{\nabla_h\Phi_M^n}
\le \mynormb{P_M\Phi^n}_{H^1}\le c_2\quad\text{for $1\le n\le N$,}
\end{align}
where $c_2$ is dependent on the domain $\Omega$
and initial data $\Phi(t_0)$, but independent of the time $t_n$.

We are in the position to prove the $L^2$ norm convergence of the adaptive BDF2 scheme \eqref{scheme: BDF2 convex CH}.
In this main theorem, $c_3:=c_{\Omega}^2(c_1^2+c_1c_2+c_2^2)$, $c_4:=16\kappa c_z^4c_3^4\mathfrak{m}_2^3\mathfrak{m}_3^2/(\mathfrak{m}_1^5\epsilon^6)$,
$c_5:=288\kappa/\bra{\mathfrak{m}_1\epsilon^2}$ and $c_{\epsilon}:=2 \bra{c_4+c_5}$.
These fixed constant may be dependent on the given data, the solution and
the starting values, but are always independent of the time $t_n$, time-step sizes $\tau_n$
and step ratios $r_n$. Moreover, they remain bounded even
when $r_n$ approach the user limit $r_{\mathrm{user}}$.

\begin{theorem}\label{thm: L2 Convergence-CH}
Assume that the CH problem \eqref{cont: Problem-CH} has a smooth solution
$\Phi\in C^3\brab{[0,T];{H}_{per}^{m+4}}$
for some integer $m\ge0$. Suppose further that the step-ratios condition $\mathbf{S0}$ and
the stabilized constraint \eqref{ieq: stabilized A restriction} hold such that
the convex-splitting BDF2 scheme \eqref{scheme: BDF2 convex CH} is unique solvable and energy stable.
If $\tau\le 1/c_{\epsilon}$,
the solution $\phi^n$ is robustly convergent in the $L^2$ norm,
\begin{align*}
\mynormb{\Phi^n-\phi^n}
\le C_{\phi}&\,\exp\brab{c_{\epsilon}t_{n-1}}\bigg(
\mynormb{\Phi_M^1-\phi^1}
+\tau\mynormb{\partial_\tau(\Phi_M^1-\phi^1)}
+t_nh^m\\
&\,+t_n\tau^2\max_{0<t\le T}\bra{
\mynormb{\Phi(t)}_{H^4}+\mynormb{\Phi''(t)}_{H^2}
+\mynormb{\Phi'''(t)}_{L^2}}\bigg)\quad\text{for $2\le n\le N$}.
\end{align*}
\end{theorem}

\begin{proof}
We evaluate the $L^2$ norm error $\mynorm{\Phi^n-\phi^n}$ by a usual splitting,
$$\Phi^n-\phi^n=\Phi^n-\Phi_M^n+e^n,$$
where $e^n:=\Phi_M^n-\phi^n\in \mathbb{\mathring V}_{h}$
is the difference between the projection $\Phi_M^n$
and the numerical solution $\phi^n$ of the convex-splitting BDF2 scheme
\eqref{scheme: BDF2 convex CH}.
Actually,  the projection solution $\Phi_M^n\in\mathscr{F}_M$,
the volume conservative  property becomes available at the discrete level
\[
\myinnerb{\Phi_M^n,1}
=\myinnerb{\Phi_M^0,1}=\myinnerb{\phi^0,1}
=\myinnerb{\phi^n,1},
\]
which implies the error function $e^n\in \mathbb{\mathring V}_{h}$.
Applying Lemma \ref{lem:Projection-Estimate}, one has
$$\mynorm{\Phi^n-\Phi_M^n}=\mynorm{I_M\brat{\Phi^n-\Phi_M^n}}_{L^2}
\le C_{\phi}\mynorm{I_M\Phi^n-\Phi_M^n}_{L^2}\le C_{\phi}h^{m}\absb{\Phi^n}_{H^{m}}.$$
Once an upper bound of $\mynorm{e^n}$ is available, the claimed error estimate follows immediately,
\begin{align}\label{Triangle-Projection-Estimate}
\mynorm{\Phi^n-\phi^n}\le\mynorm{\Phi^n-\Phi_M^n}+\mynorm{e^n}
\le  C_{\phi}h^{m}\absb{\Phi^n}_{H^{m}}+\mynorm{e^n}
\quad \text{for $1\le n\le N$}.
\end{align}

To bound $\mynorm{e^n}$, we consider two stages:
Stage 1 analyzes the space consistency error for a semi-discrete system having a projected solution $\Phi_M$;
With the help of the Young-type and embedding convolution inequalities with respect to DOC kernels $\theta_{k-j}^{(k)}$
and the solution estimate in Lemma \ref{lem: bound-BDF2 Solution CH},
Stage 2 derives the error estimate for the fully discrete error system.

\paragraph{Stage 1: Consistency analysis of semi-discrete projection}
A substitution of the projection solution $\Phi_M$ and differentiation  operator $\Delta_h$ into the original equation \eqref{cont: Problem-CH} yields the semi-discrete system
\begin{align}\label{Projection-Equation}
\partial_t\Phi_M
=\kappa\Delta_h\mu_M
+\zeta_{P}\quad\text{with}\quad
\mu_M=F^\prime(\Phi_M)-\epsilon^2\Delta_h\Phi_M,
\end{align}
where $\zeta_{P}(\myvec{x}_h,t)$ represents the spatial consistency error arising from the projection of exact solution, that is,
\begin{align}\label{Projection-truncation error}
\zeta_{P}:=\partial_t\Phi_M-\partial_t\Phi
+\kappa(\Delta\mu - \Delta_h\mu_M)\quad \text{for $\myvec{x}_h\in\Omega_{h}$.}
\end{align}

Following the proof of \cite[Theorem 3.1]{LiaoJiZhang:2020pfc},
and using Lemma \ref{lem:Projection-Estimate}, it is not difficult to obtain that $\mynorm{\zeta_{P}}\le C_\phi h^m$
and $\mynorm{\zeta_{P}(t_j)}\le C_\phi h^m$ for $j\ge2$.
Then Lemma \ref{lem: DOC property} (III) yields
\begin{align}\label{Projection-consistency}
\sum_{k=2}^n\mynormb{\Upsilon_{P}^k}\le C_\phi h^m \sum_{k=2}^n\sum_{j=2}^k\theta_{k-j}^{(k)}
\le C_\phi t_nh^m\quad\text{where}\;\; \Upsilon_{P}^k:=\sum_{j=2}^k\theta_{k-j}^{(k)}\zeta_{P}(t_j)\;\;\text{for $k\ge2$.}
\end{align}

\paragraph{Stage 2: $L^2$ norm error of fully discrete system}
From the projection equation \eqref{Projection-Equation},
one can apply the BDF2 formula to obtain the following approximation equation
\begin{align}\label{Discrete-Projection-Equation}
D_2\Phi_M^n
=\kappa\Delta_h\mu_M^n +\zeta_{P}^n + \xi_{\Phi}^n\quad\text{with}\quad
\mu_M^n=
\brab{\Phi_M^n}^3-\hat{\Phi}_M^n-\bra{\epsilon^2+A\tau^2}\Delta_h\Phi_M^n,
\end{align}
where the local consistency errors $\xi_{\Phi}^n$
and $\zeta_{P}^n:=\zeta_{P}(t_n)$ are defined by
\eqref{def: BDF2-local consistency} and \eqref{Projection-truncation error},
respectively.
Subtracting the full discrete scheme \eqref{scheme: BDF2 convex CH} from the approximation equation \eqref{Discrete-Projection-Equation},
we have the following error system
\begin{align}\label{Error-Equation}
D_2e^n
=\kappa\Delta_h\kbrab{f_{\phi}^ne^n-\hat{e}^n-\bra{\epsilon^2+A\tau^2}\Delta_he^n}
+\zeta_{P}^n+\xi_{\Phi}^n\quad\text{for $2\le n\le N$,}
\end{align}
where the nonlinear term $f_{\phi}^n :=\brat{\Phi_M^n}^2+\Phi_M^n\phi^n+\brat{\phi^n}^2$
and $\hat{e}^n:=e^{n-1}-e^{n-2}$.
Thanks to the estimates in Lemma \ref{lem: bound-BDF2 Solution CH} and \eqref{H1 bound projection},
one applies the embedding inequality \eqref{ieq: H1 embedding L6} to find that
\begin{align}\label{ieq: L3 bound nonlinear}
\mynormb{f_{\phi}^n}_{l^3}\le&\, \mynormb{\Phi_M^n}_{l^6}^2
+\mynormb{\Phi_M^n}_{l^6}\mynormb{\phi^n}_{l^6}
+\mynormb{\phi^n}_{l^6}^2\le c_3.
\end{align}

Multiplying both sides of equation \eqref{Error-Equation} by the DOC kernels $\theta_{k-n}^{(k)}$,
and summing up $n$ from $n=2$ to $k$, we apply the equality \eqref{eq: orthogonal equality for BDF2} with $v^j=e^j$ to obtain
\begin{align}\label{Error-Equation-DOC}
\diff e^k
=-\theta_{k-2}^{(k)}b_{1}^{(2)}\diff e^1+\kappa\sum_{j=2}^k\theta_{k-j}^{(k)}\Delta_h
\kbrab{f_{\phi}^j e^j-\hat{e}^j-\bra{\epsilon^2+A\tau^2}\Delta_he^j}
+\Upsilon_{P}^k+\Xi_{\Phi}^k
\end{align}
for $2\le k\le N$, where $\Xi_{\Phi}^k$ and $\Upsilon_{P}^k$ are
defined by \eqref{def: BDF2-global consistency} and \eqref{Projection-consistency}, respectively.
Making the inner product of \eqref{Error-Equation-DOC} with $2e^k$,
and summing $k$ from 2 to $n$, we obtain
\begin{align}\label{Error-Equation-Inner}
\mynormb{e^n}^2\le\mynormb{e^1}^2
-2\sum_{k=2}^n\theta_{k-2}^{(k)}b_{1}^{(2)}\mynormb{e^k}\mynormb{\diff e^1}
+J^n +2\sum_{k=2}^n\myinnerb{\Upsilon_{P}^k+\Xi_{\Phi}^k,e^k}
\end{align}
for $2\le n\le N$, where $J^n$ is defined by
\begin{align}\label{Error-quadratic forms}
J^n:=&\,2\kappa\sum_{k,j}^{n,k}\theta_{k-j}^{(k)}
\myinnerb{f_{\phi}^je^j-\hat{e}^j
-\bra{\epsilon^2+A\tau^2}\Delta_he^j,\Delta_h e^k}.
\end{align}
Taking $u^j:=f_{\phi}^j$ (with the upper bound $c_u:=c_3$),
$w^j:=e^j$ and $\varepsilon=\epsilon^2/2$
in Lemma \ref{lem:DOC quadr form H2-H2 embedding inequ},
one applies the solution bound \eqref{ieq: L3 bound nonlinear} to obtain
\begin{align*}
2\kappa\sum_{k,j}^{n,k}\theta_{k-j}^{(k)} \myinnerb{f_{\phi}^je^j,\Delta_h e^k}
\le
\kappa\epsilon^2\sum_{k,j}^{n,k}\theta_{k-j}^{(k)}\myinnerb{\Delta_h e^j,\Delta_h e^k}
+c_4\sum_{k=2}^n\tau_k  \mynormb{e^k}^2.
\end{align*}
For the second term of \eqref{Error-quadratic forms},
one applies the Young-type convolution inequality
in Lemma \ref{lem:DOC quadr form Young inequ2}
by taking $w^k:=\Delta_h e^k$, $v^j:=-\hat{e}^j$
and $\varepsilon=\epsilon^2/2$ to get
\begin{align*}
2\kappa\sum_{k,j}^{n,k}\theta_{k-j}^{(k)} \myinnerb{-\hat{e}^j,\Delta_h e^k}
\le \kappa\epsilon^2\sum_{k,j}^{n,k}\theta_{k-j}^{(k)}\myinnerb{\Delta_h e^j,\Delta_h e^k}
+c_5\sum_{k=1}^{n-1}\tau_k\mynormb{e^k}^2.
\end{align*}
An application of the positive definiteness of the
kernels $\theta_{k-j}^{(k)}$ in Lemma \ref{lem: DOC property} (I) yields
\begin{align*}
2\kappa\sum_{k,j}^{n,k}\theta_{k-j}^{(k)}
\myinnerb{-A\tau^2\Delta_he^j,\Delta_h e^k}<0.
\end{align*}
Then the term $J^n$ in \eqref{Error-quadratic forms} can be bounded by
\begin{align*}
J^n\le&\,\frac{c_{\epsilon}}{2}\sum_{k=1}^n\tau_k  \mynormb{e^k}^2.
\end{align*}

Therefore, it follows from \eqref{Error-Equation-Inner} that
\begin{align*}
\mynormb{e^n}^2
\le \mynormb{e^1}^2-2\sum_{k=2}^n\theta_{k-2}^{(k)}b_{1}^{(2)}\mynormb{e^k}\mynormb{\diff e^1}
+\frac{c_{\epsilon}}{2}\sum_{k=1}^n\tau_k  \mynormb{e^k}^2
+ 2\sum_{k=2}^n\mynormb{e^k}\mynormb{\Upsilon_{P}^k+\Xi_{\Phi}^k}
\end{align*}
for $2\le n\le N$. Choosing some integer $n_0$ ($1\le n_0 \le n$) such that
$\mynormb{e^{n_0}}=\max_{1\le k \le n}\mynormb{e^k}$.
Taking $n:=n_0$ in the above inequality, one can obtain
\begin{align*}
\mynormb{e^{n_0}}\le \mynormb{e^1}-2\mynormb{\partial_{\tau} e^1}\sum_{k=2}^{n_0}\theta_{k-2}^{(k)}b_{1}^{(2)}\tau_1
+\frac{c_{\epsilon}}{2}\sum_{k=1}^{n_0}\tau_k  \mynormb{e^k}+ 2\sum_{k=2}^{n_0}
\mynormb{\Upsilon_{P}^k+\Xi_{\Phi}^k}.
\end{align*}
By using Lemma \ref{lem: DOC property} (II), one has
$$-\theta_{k-2}^{(k)}b_{1}^{(2)}\tau_1=\tau_1\prod_{i=2}^k\frac{r_i^2}{1+2r_i}=\tau_k\prod_{i=2}^k\frac{r_i}{1+2r_i}
\leq \frac{\tau_k}{2^{k-1}}\quad \text{for $2\le k\le N$},$$
such that
$$-\sum_{k=2}^{n}\theta_{k-2}^{(k)}b_{1}^{(2)}\tau_1
\le\tau\sum_{k=2}^{n}\frac{1}{2^{k-1}}\le \tau\quad\text{for $2\le n\le N$}.$$
Thus one gets
\begin{align*}
\mynormb{e^n}\le\mynormb{e^{n_0}}
&\le \mynormb{e^1}+2\tau\mynormb{\partial_{\tau} e^1}
+\frac{c_{\epsilon}}{2}\sum_{k=1}^{n}\tau_k  \mynormb{e^k}
+ 2\sum_{k=2}^{n}\mynormb{\Upsilon_{P}^k+\Xi_{\Phi}^k}.
\end{align*}
Under the maximum step constraint $\tau\le 1/c_{\epsilon}$, we have
\begin{align*}
\mynormb{e^n}\le 2(1+c_{\epsilon})\mynormb{e^1}
+4\tau\mynormb{\partial_{\tau} e^1}
+c_{\epsilon}\sum_{k=2}^{n-1}\tau_k  \mynormb{e^k}
+ 4\sum_{k=2}^n\mynormb{\Upsilon_{P}^k+\Xi_{\Phi}^k}.
\end{align*}
The discrete Gr\"onwall inequality \cite[Lemma 3.1]{LiaoZhang:2020linear} yields the following estimate
\begin{align*}
\mynormb{e^n}&\le 2\exp\brab{c_{\epsilon} t_{n-1}}
\kbraB{(1+c_{\epsilon})\mynormb{e^1}
+2\tau\mynormb{\partial_{\tau} e^1}
+2\sum_{k=2}^n\mynormb{\Upsilon_{P}^k}+ 2\sum_{k=2}^n\mynormb{\Xi_{\Phi}^k}}
\end{align*}
for $2\le n\le N$.
Furthermore, the convolutional consistency error established
in Lemma \ref{lem: BDF2-Consistency-Error} together with
the regularity condition $\Phi\in C^3\brab{[0,T];{H}_{per}^{m+4}}$
and Lemma \ref{lem:Projection-Estimate},
gives the bound of the global temporal error term $\sum_{k=2}^n\mynormb{\Xi_{\Phi}^k}$.
Therefore by applying the error estimate \eqref{Projection-consistency}
and the triangle inequality \eqref{Triangle-Projection-Estimate}, we complete the proof.
\end{proof}

\section{Numerical experiments}
\setcounter{equation}{0}
We run the BDF2 scheme \eqref{scheme: BDF2 convex CH}
for the CH equation \eqref{cont: Problem-CH}.
In our computations,
the parameter $A=3/625$ according to
Remark \ref{remark: comments on stabilized A condition}.
The TR-BDF2 method is always employed to obtain the first-level solution.
A simple fixed-point iteration with the termination error $10^{-12}$
is employed to solve the nonlinear algebra equations at each time level.

\subsection{Robustness tests on random time meshes}
\begin{example}\label{exam:CH BDF2 Accuracy Test}
To facilitate the robustness test of the convex-splitting BDF2 method
\eqref{scheme: BDF2 convex CH},
we consider an exact solution
$\Phi(\mathbf{x},t)=\cos(t)\sin({x})\sin({y})$
with the model parameters $\kappa=2\times10^{-3}$ and $\epsilon=5\times10^{-2}$
by adding a corresponding exterior force to the CH model \eqref{cont: Problem-CH}.
\end{example}

In the following  examinations,
the computational domain $(0,2\pi)^2$ is discretized
by using $128^2$ spatial meshes.
Then the problem is solved until time $T=1$ on random time meshes.
To be more precise,
we take the time step sizes $\tau_k:=T\sigma_{k}/S$ for $1\leq k\leq N$,
where $\sigma_{k}\in(0,1)$ is the uniformly distributed random number
and $S=\sum_{k=1}^N\sigma_{k}$.
Since the spectral accuracy in space is standard, we only test the time accuracy
with the numerical error $e(N):=\max_{1\leq{n}\leq{N}}\mynorm{\Phi(t_n)-\phi^n}$
in each run.  The numerical order of convergence is estimated by
$\text{Order}:=\log\bra{e(N)/e(2N)}/{\log\bra{\tau(N)/\tau(2N)}},$
where $\tau(N)$ denotes the maximum time-step size for total $N$ subintervals.

\begin{table}[htb!]
\begin{center}
\caption{Accuracy of BDF2 method \eqref{scheme: BDF2 convex CH} on random time meshes.}\label{examp:CH BDF2 L2 error} \vspace*{0.3pt}
\def\temptablewidth{0.7\textwidth}
{\rule{\temptablewidth}{0.5pt}}
\begin{tabular*}{\temptablewidth}{@{\extracolsep{\fill}}cccccc}
  $N$   &$\tau$      &$e(N)$     &Order  &$\max r_k$ &$N_1$\\
  \midrule
  40	 &3.96e-02	 &3.69e-04	 &1.94	 &17.27	     &3\\	
  80	 &2.44e-02	 &1.08e-04	 &2.55	 &46.22	     &5\\	
  160	 &1.29e-02	 &2.75e-05	 &2.13	 &167.41	 &16\\	
  320	 &6.28e-03	 &7.07e-06	 &1.90	 &264.04	 &29\\	
  640	 &3.05e-03	 &1.57e-06	 &2.08	 &1584.01	 &62\\			
\end{tabular*}
{\rule{\temptablewidth}{0.5pt}}
\end{center}
\end{table}

The numerical results obtained using a set of random meshes
are tabulated in Table  \ref{examp:CH BDF2 L2 error}.
In addition to the discrete $L^2$ numerical error between the exact solution
and the numerical solution,
the maximum  time-step size $\tau$,
the maximum step ratio $\max r_k$ and the number (denote by $N_1$)
of time levels with the step ratios $r_k\ge 4.864$
are also recorded, respectively.

As observed, the convex-splitting BDF2 method
\eqref{scheme: BDF2 convex CH} still achieves
the second-order accuracy on arbitrary nonuniform meshes even though
some step ratios lager than  $r_{*}\approx 4.864$.
The numerical results indicate that
the BDF2 method is robust with respect to
the step-size variations than previous theoretical predictions.
Also, the improved condition $0<r_k<4.864$ is still a
sufficient condition for second-order convergence.

\begin{figure}[htb!]
\centering
\subfigure[time-step size $\tau=10^{-1}$]{
\includegraphics[width=2.0in]{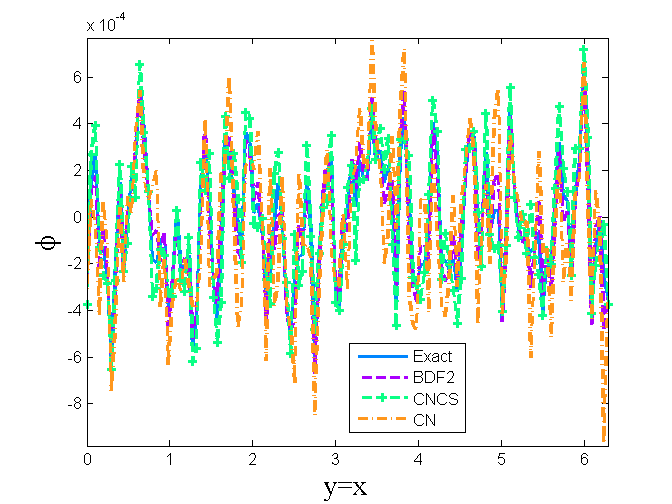}}
\subfigure[time-step size $\tau=5\times10^{-2}$]{
\includegraphics[width=2.0in]{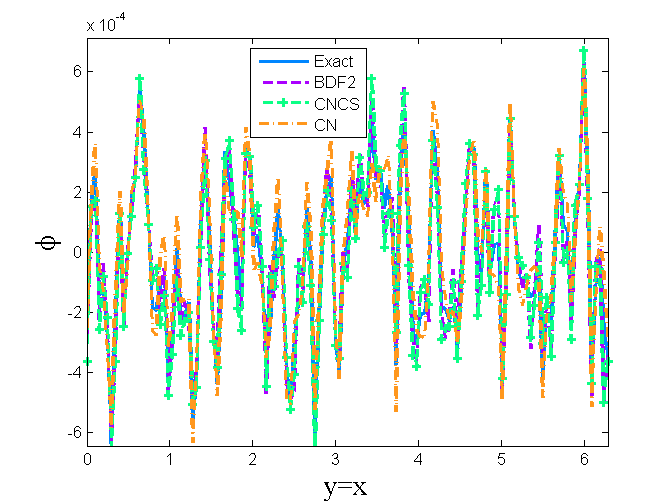}}
\subfigure[time-step size $\tau=10^{-2}$]{
\includegraphics[width=2.0in]{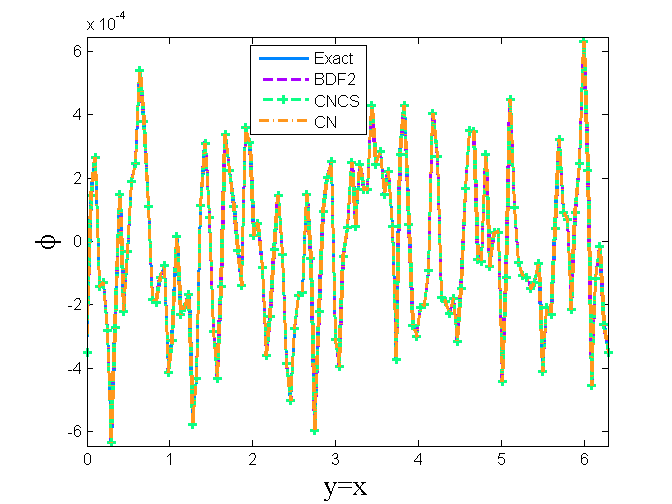}}\\
\caption{Solution curves by BDF2, CN and CNCS methods at $T=0.1$.}
\label{examp:compar solution}
\end{figure}
\begin{figure}[htb!]
\centering
\subfigure[time-step size $\tau=10^{-1}$]{
\includegraphics[width=2.0in]{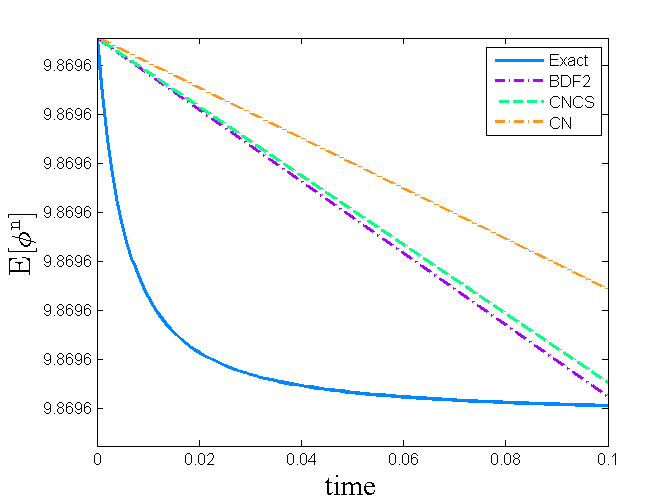}}
\subfigure[time-step size $\tau=5\times10^{-2}$]{
\includegraphics[width=2.0in]{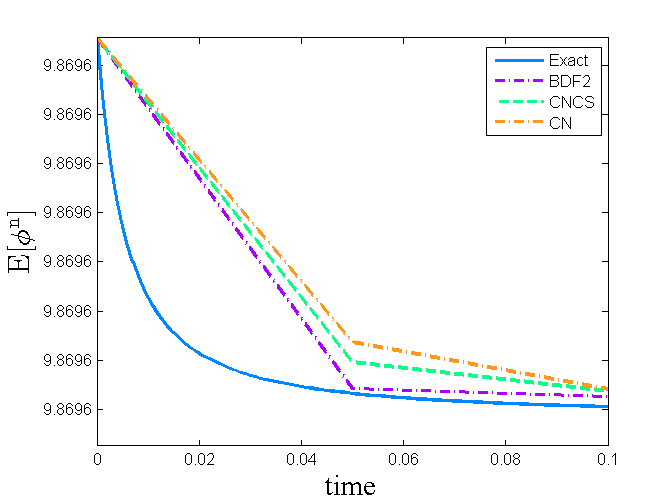}}
\subfigure[time-step size $\tau=10^{-2}$]{
\includegraphics[width=2.0in]{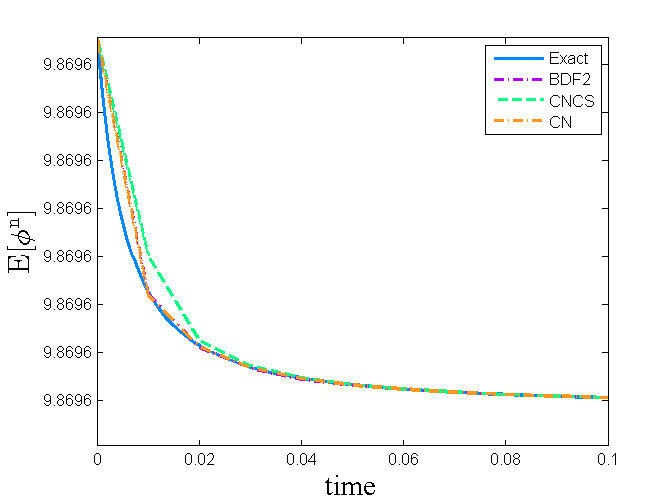}}\\
\caption{Original energy curves by BDF2, CN and CNCS methods until $T=0.1$.}
\label{examp:compar energy}
\end{figure}

\begin{example}\label{exam:coarsen dynamic}
We next simulate the coarsening dynamics of the
CH equation \eqref{cont: Problem-CH}.
Precisely, the initial condition is taken as
$\Phi_0(\mathbf{x})=\text{rand}(\mathbf{x})$,
where $\text{rand}(\mathbf{x})$ generates random numbers
between $-0.001$ to $0.001$ uniformly.
Here, the mobility coefficient $\kappa=2\times10^{-3}$ and
the interfacial thickness $\epsilon=5\times10^{-2}$ are taken
in the following numerical simulations.
Always, the spatial domain $(0,2\pi)^2$ is discretized
by using $128^2$ spatial meshes.
\end{example}

\subsection{Numerical comparisons}
To further benchmark the convex-splitting BDF2 scheme with the random initial data
generated in Example \ref{exam:coarsen dynamic},
we run several numerical tests to explore the numerical behaviors
near the initial time.
We also implement the unconditionally energy stable
Crank-Nicolson (CN) method \cite{ZhangQiao:2012},
\begin{align*}
\partial_\tau\phi^n
&=\kappa\Delta_h\mu^{n-\frac12}\quad\text{with}\quad
\mu^{n-\frac12}=\frac{1}{2}\kbrab{(\phi^{n})^2
+(\phi^{n-1})^2}\phi^{n-\frac12}
-\phi^{n-\frac12}-\varepsilon^2\Delta_h\phi^{n-\frac12},
\end{align*}
and the second-order Crank-Nicolson convex-splitting
(CNCS) method \cite{ChengWangWiseYue:2016Weakly,GuoWangWiseYue:2016CHSecondConvex},
\begin{align*}
\partial_\tau\phi^n
&=\kappa\Delta_h\hat{\mu}^{n-\frac12}\quad\text{with}\quad
\hat{\mu}^{n-\frac12}
=\frac{1}{2}\kbrab{(\phi^{n})^2+(\phi^{n-1})^2}\phi^{n-\frac12}
-\check{\phi}^{n-\frac12}
-\varepsilon^2\Delta_h\hat{\phi}^{n-\frac12},
\end{align*}
where $\phi^{n-\frac12}:=(\phi^n+\phi^{n-1})/2$, $\hat{\phi}^{n-\frac12}:=\bra{3\phi^n+\phi^{n-2}}/4$
and $\check{\phi}^{n-\frac12}=\bra{3\phi^{n-1}-\phi^{n-2}}/2$.
Since the CNCS method requires two initialization steps,
a first-order convex-splitting
scheme \cite{ChenWangYanZhang:2019} is used here to obtain the first-level solution.

The random initial data initiates a fast coarsening dynamics at the beginning time.
We use a random initial profile to test the effectiveness of various numerical methods
with different time step sizes.
The numerical solution curves are summarized in Figure  \ref{examp:compar solution},
where the reference solution
is obtained by using the convex-splitting BDF2 method
with a uniform time-step size $\tau=10^{-3}$.
We observe that solutions of CN and CNCS methods tend to
generate non-physical oscillations when some large time steps are used.
In contrast, the convex-splitting BDF2 solution
is more robust and accurate than the CN and CNCS schemes with the same time step size.
It seems that
the BDF2 method is more suitable than Crank-Nicolson type schemes
when large time-step sizes are adopted.

\begin{figure}[htb!]
\centering
\subfigure[Original energy]{
\includegraphics[width=2.7in,height=1.8in]{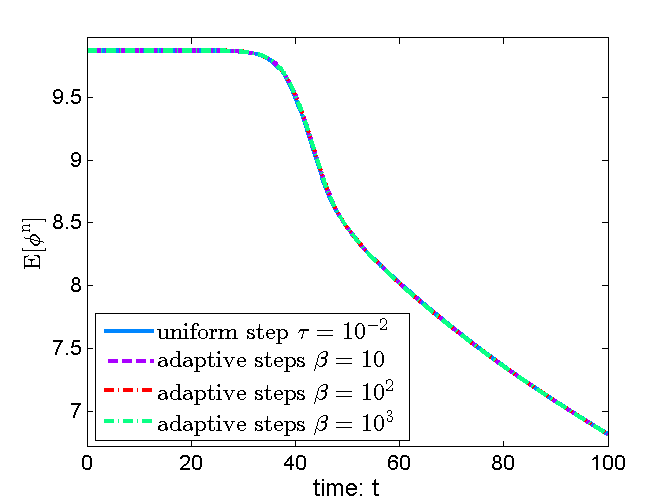}}
\subfigure[Time step sizes]{
\includegraphics[width=2.7in,height=1.8in]{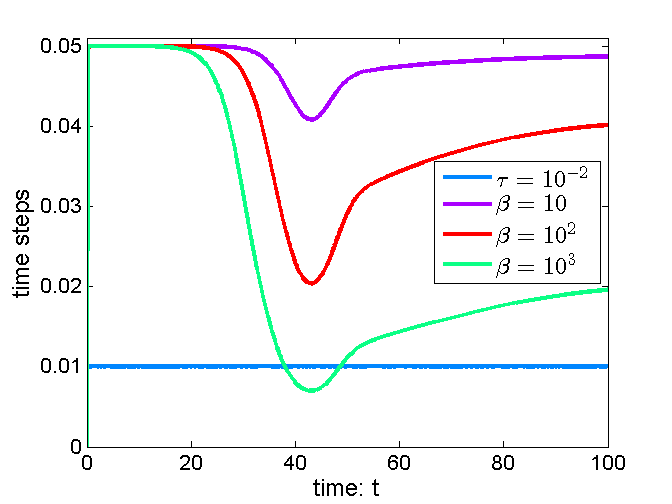}}
\caption{Energy curves and adaptive time-step sizes for different parameters $\beta$.}
\label{examp:compar adap param}
\end{figure}

\begin{table}[htb!]
\begin{center}
\caption{CPU time (in seconds) and total time steps comparisons.}
\label{examp:compar adap param CPU} \vspace*{0.3pt}
\def\temptablewidth{0.6\textwidth}
{\rule{\temptablewidth}{0.5pt}}
\begin{tabular*}{\temptablewidth}{@{\extracolsep{\fill}}c|c|c|c|c}
  Strategies   &$\tau=10^{-2}$ &$\beta=10$  &$\beta=10^2$  & $\beta=10^3$ \\
  \midrule
  CPU time    &109.116   &35.601     &39.238       &71.880\\
  Time levels &10000     &2098       &2710         &5671\\
\end{tabular*}
{\rule{\temptablewidth}{0.5pt}}
\end{center}
\end{table}	

\subsection{Simulation of coarsening dynamics}
In this subsection, we simulate the coarsening dynamics by using
the convex-splitting BDF2 method \eqref{scheme: BDF2 convex CH}
with the random initial condition. In what follows,
to capture the multiple time scales accurately
and to improve the computational efficiency for long-time simulations,
the time steps are selected by using
the following adaptive time-stepping strategy
\cite{HuangYangWei:2020},
\begin{align}\label{algo:adaptive step}
\tau_{ada}
=\max\Bigg\{\tau_{\min},
\frac{\tau_{\max}}{\sqrt{1+\beta\mynormb{\partial_\tau \phi^n}^2}}\Bigg\}
\quad\text{so that}\quad
\tau_{n+1}=\min\big\{\tau_{ada},r_{\mathrm{user}}\tau_n\big\},
\end{align}
where $\beta>0$ is a user chosen parameter,
$\tau_{\max}$ and $\tau_{\min}$ are
the predetermined maximum and minimum time steps, respectively.

We take $r_{\mathrm{user}}=4$, $\tau_{\min}=5\times10^{-5}$ and $\tau_{\max}=5\times10^{-2}$ in the adaptive time-stepping algorithm \eqref{algo:adaptive step},
and run the convex-splitting BDF2 method
\eqref{scheme: BDF2 convex CH} until time $T=100$.
The reference solution is obtained by applying a small time step $\tau=10^{-2}$.
As seen in Figure \ref{examp:compar adap param}, we use three different user parameters
$\beta=10,10^2$ and $10^3$ to compute the discrete original energy and the corresponding adaptive time-steps.
One can observe that the discrete energy curves
using the adaptive stepping algorithm are comparable to the reference one.
On the other hand, the adjustments of time-steps are closely relied on
the user parameter $\beta$. As expected, a large $\beta$ leads to small time-step sizes,
and a small $\beta$ generates large step sizes.
The CPU time (in seconds) and the adaptive time levels
recorded in Table \ref{examp:compar adap param CPU}
show the effectiveness and efficiency of the adaptive time-stepping algorithm,
which makes the long-time dynamics simulations practical.

\begin{figure}[htb!]
\centering
\subfigure[time $t=10$]{
\includegraphics[width=2in]{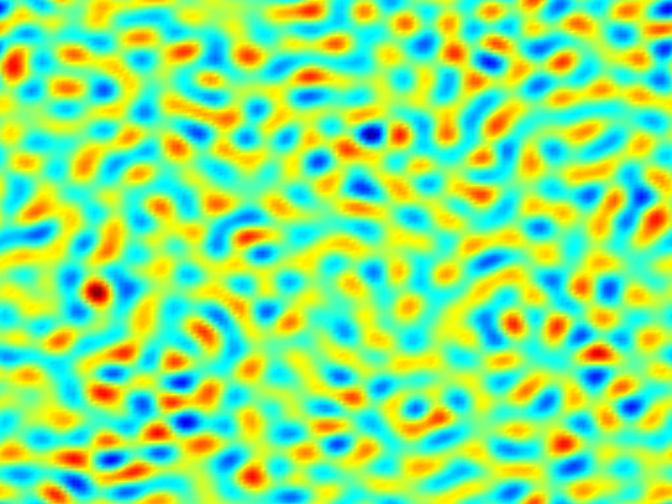}}
\subfigure[time $t=100$]{
\includegraphics[width=2in]{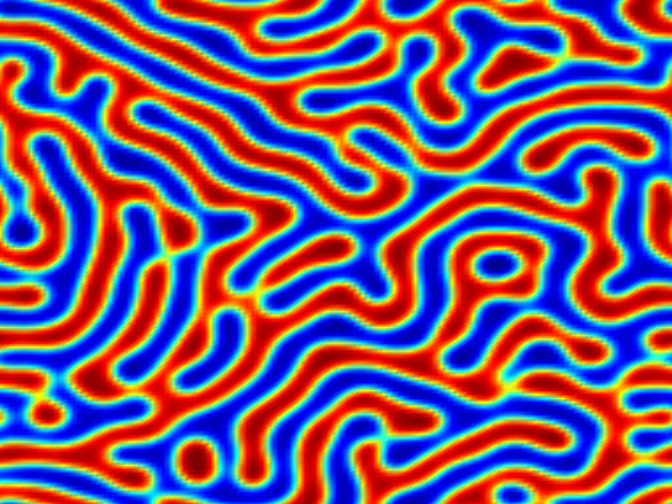}}
\subfigure[time $t=300$]{
\includegraphics[width=2in]{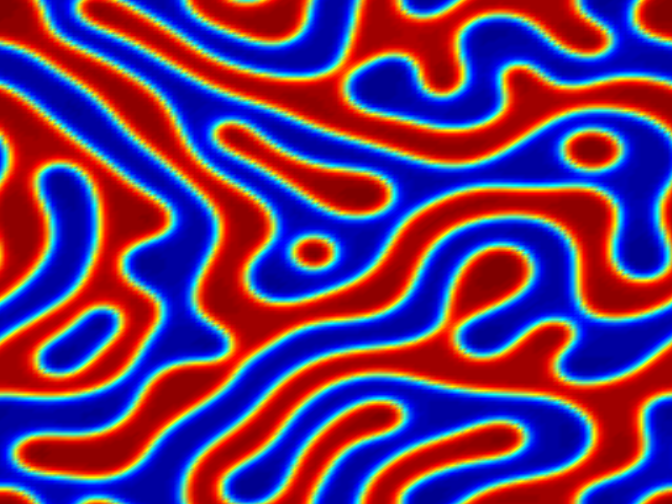}}\\
\subfigure[time $t=500$]{
\includegraphics[width=2in]{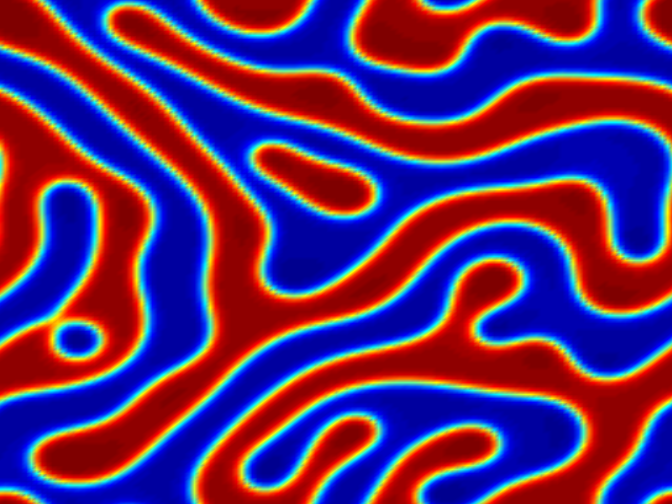}}
\subfigure[time $t=800$]{
\includegraphics[width=2in]{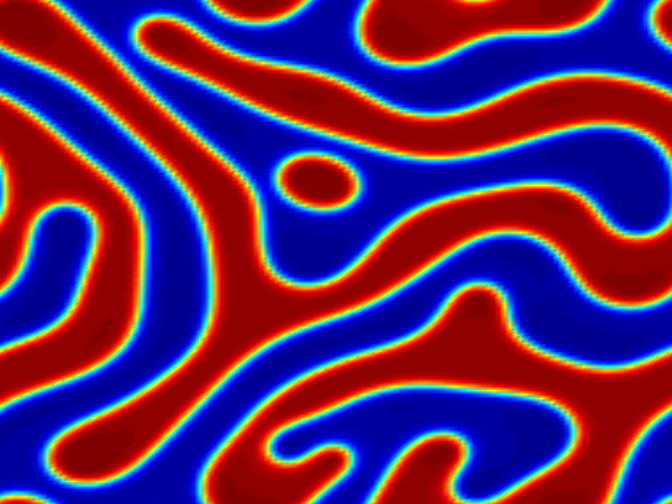}}
\subfigure[time $t=1000$]{
\includegraphics[width=2in]{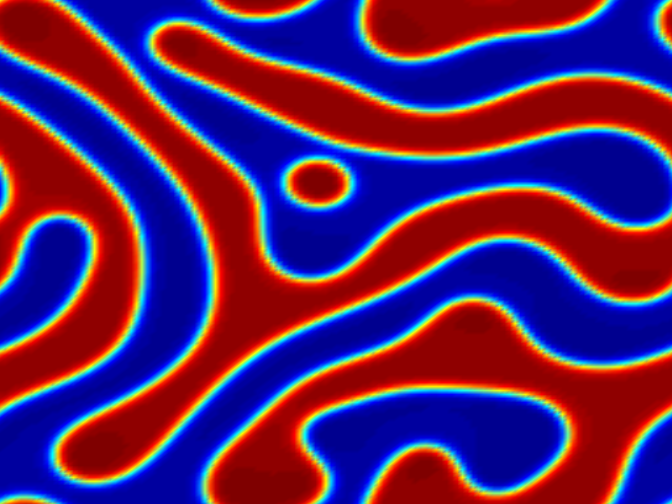}}\\
\caption{The profile of numerical solution $\phi$ at different time for the CH model.}
\label{examp:CH coarsen snap}
\end{figure}

\begin{figure}[htb!]
\centering
\subfigure[Energy scaling]{
\includegraphics[width=2.0in]{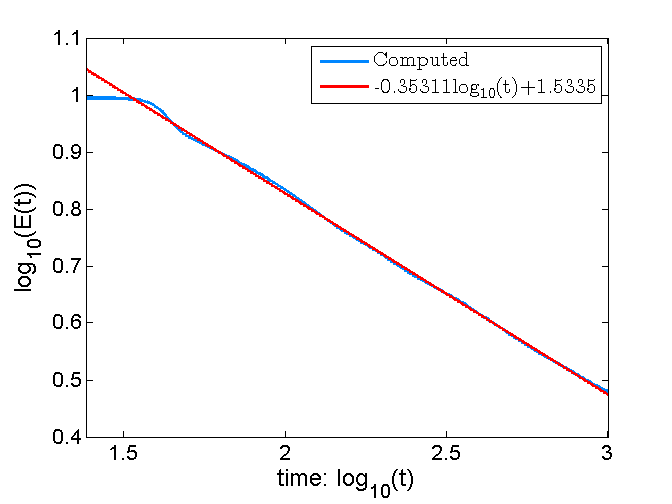}}
\subfigure[Volume difference]{
\includegraphics[width=2.0in]{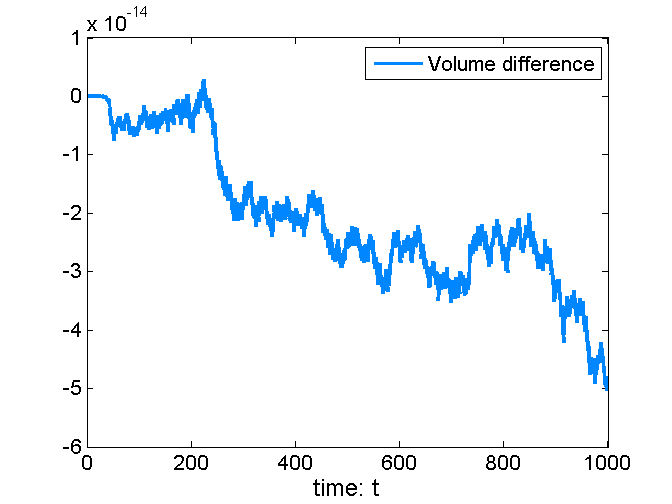}}
\subfigure[Adaptive step sizes]{
\includegraphics[width=2.0in]{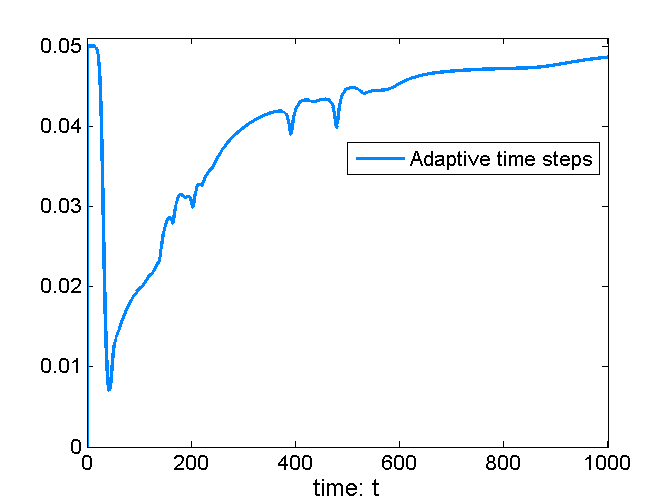}}\\
\caption{Numerical results show original energy, volume and adaptive
  time steps of the CH equation during the coarsening dynamics.}
\label{examp:adapt energy volume}
\end{figure}

We next perform the coarsening dynamic simulations by using the
above adaptive time-stepping strategy with the setting $\beta=10^3$
until time $T=1000$.
The evolution of microstructure for the CH
model due to the phase separation  at different time
are summarized in Figure \ref{examp:CH coarsen snap}.
As seen, the microstructure is relatively fine and consists of many
precipitations at early time.
The coarsening, dissolution, merging processes are also observed.
The time evolutions of original energy, volume and the adaptive step sizes are summarized in
Figure \ref{examp:adapt energy volume}.
The subplot (a) of Figure \ref{examp:adapt energy volume}
demonstrates a very good agreement with the expected scaling law, i.e.,
the energy decreases as $O(t^{-\frac13})$.

\end{document}